\documentclass[9pt,amstex]{article}
\usepackage{amssymb,amsmath} 
\usepackage{latexsym}

\newtheorem{thm}{Theorem}[section] 
\newtheorem{dfn}[thm] {Definition}
\newtheorem{rmk}[thm]{Remark}
\newtheorem{cor}[thm]{Corollary}
\newtheorem{prop}[thm]{Proposition}

\newtheorem{lemma}[thm]{Lemma}

\newtheorem{ex}[thm]{Example}

\def\cyclic{\mathop{\kern0.9ex{{+}
\kern-2.2ex\raise-.28ex\hbox{\Large\hbox
{$\circlearrowright$}}}}\limits}

\def\buildrel#1_#2^#3{\mathrel{\mathop{\kern 0pt#1}\limits_{#2}^{#3}}}

\newcommand{\Pf}{\noindent{\em Proof.} }
\newcommand{\EPf}
{%
\mbox{}%
\nolinebreak%
\hfill%
\rule{2mm}{2mm}%
\medbreak%
\par%
}

\newcommand{\Svar}{\mbox{\rm Svar}}
\newcommand{\Ker}{\mbox{\rm Kern}}

\newcommand{\Tr}{\mbox{$\mathtt{Tr}$}}

\newcommand{\End}{\mbox{$\mathtt{End}$}}
\newcommand{\Hom}{\mbox{$\mathtt{Hom}$}}

\newcommand{\im}{\mbox{$\mathtt{Im}$}}

\newcommand{\C}{\mathbb C}

\newcommand{\A}{\mathbb A}

\newcommand{\Z}{\mathbb Z} 

\newcommand{\R}{\mathbb R}
\newcommand{\K}{\mathbb K}
 
\renewcommand{\H}{\mathbb H} 

\newcommand{\e}{{\mathfrak{e}}}

\newcommand{\g}{{\mathfrak{g}}{}}

\newcommand{\p}{{\mathfrak{p}}{}} 
\newcommand{\Der}{{\mathfrak{Der}}{}} 
\newcommand{\q}{{\mathfrak{q}}{}}

\newcommand{\h}{{\mathfrak{h}}{}}

\newcommand{\sL}{{\mathfrak{sl}}{}}

\def\cref#1{Corollary~\ref{#1}}

\newcommand{\sP}{\mathfrak{sp}}
\newcommand{\oo}{\mathfrak{o}}
\newcommand{\so}{\mathfrak{so}}
\newcommand{\gL}{\mathfrak{gl}}

\newcommand{\su}{\mathfrak{su}}

\newcommand{\Gl}{\mathrm{GL}}
\newcommand{\Sl}{\mathrm{SL}}
\newcommand{\UU}{\mathrm{U}}

\newcommand{\OO}{\mathrm{O}}
\newcommand{\SO}{\mathrm{SO}}
\newcommand{\Sp}{\mathrm{Sp}}
\newcommand{\Herm}{\mathrm{Herm}}
\newcommand{\Sym}{\mathrm{Sym}}

\newcommand{\Asym}{\mathrm{Asym}}

\newcommand{\id}{\mathrm{id}}

\newcommand{\cX}{\mathcal{X}}

\newcommand{\HHH}{\mathbb H} 
\newcommand{\OOO}{\mathbb{O}}
\newcommand{\PPP}{\mathbb{P}}
\newcommand{\bF}{\mathbb{F}}
\newcommand{\msk}{\medskip}
\newcommand{\ssk}{\smallskip}
\newcommand{\nin}{\noindent}

\title{Homotopes of Symmetric Spaces\\
II. Structure Variety and Classification}

\author{
{\bf Wolfgang Bertram}\\
IECN Universit\'e Nancy I, France\\
e-mail:  {\tt bertram@iecn.u-nancy.fr}\\
{\bf Pierre Bieliavsky}\\
Universit\'e Catholique de Louvain, Belgium.\\
e-mail: {\tt pierre.bieliavsky@gmail.com}\\
}
\date{}%10 november  2010} 

\begin{document}
\maketitle

{\small
\nin {\bf Abstract.}
We classify {\em homotopes} of classical symmetric spaces (introduced in \cite{Be08} and
studied in Part I of this work). Our classification uses the {\em fibered structure} of homotopes:
they are fibered as symmetric spaces, with flat fibers, over a non-degenerate base; the
base spaces correspond to {\em inner ideals} in Jordan pairs. Using that inner ideals 
in classical Jordan pairs are always complemented in the sense of \cite{LoN94}, the 
classification of homotopes is obtained by combining the classification of inner ideals
with the one of isotopes of a given inner ideal.

\msk

{\bf \nin Subject classification.}
17B60,  %	Lie (super)algebras associated with other structures (associative, Jordan, etc.)
17B70  %	Graded Lie (super)algebras
17C37,  % 	Jordan: Associated geometries
53C35  %	Symmetric spaces [See also 
32M15 % Hermitian symmetric spaces, bounded symmetric domains, Jordan algebras

\msk
\nin {\bf Key words.}
contraction, complementation, involution, inner ideal, homotope (isotope), 
Jordan pairs and triple systems,
$3$-graded Lie algebra, Lie triple system, structural transformation, symmetric space,
pseudo-inverse, idempotent}

\section*{Introduction}

{\em Triple systems} are trilinear maps $V^3 \to V$ satisfying certain algebraic
properties, giving rise to notions such as
 {\em Lie-, Jordan-} or {\em associative triple systems}. 
Often it is more conceptual to consider
a pair of vector spaces $(V^+,V^-)$ and two trilinear maps
$T^\pm: V^\pm \times V^\mp \times V^\pm \to V^\pm$,
verifiying certain algebraic properties (e.g., {\em Jordan pairs}, {\em associative pairs}).
Given such an object $(V^\pm,T^\pm)$,
we define  a subset of $\Hom(V^+,V^-)$, called the {\em structure variety
of $(V^\pm,T^\pm)$}:
\begin{equation}\label{SvarDefEqn}
\Svar(V^+):=\big\{ \alpha:V^+ \to V^- | \, \alpha \mbox{ linear, }   \forall u,v,w \in V^+ :
T^- (\alpha u,v,\alpha w)=\alpha T^+(u,\alpha v,w)  \big\}.
\end{equation}
Then elements $\alpha \in \Svar(V^+)$
parametrize families of trilinear maps on $V^+$, given by
\begin{equation}\label{THomotopeEqn}
T_\alpha : V^+ \times V^+ \times V^+ \to V^+, \quad T_\alpha(x,y,z):=T^+(x,\alpha(y),z) \ .
\end{equation}
It is well-known (cf.\ Lemma \ref{StructureLemma})  that,
if $(V^\pm,T^\pm)$ is a {\em Jordan pair}, then $T_\alpha$ is a {\em Jordan triple system},
and since skew-symmetrizing a Jordan triple system yields a Lie triple system
(Lemma \ref{LieJordanLemma}), we conclude that
\begin{equation}\label{RHomotopeEqn}
R_\alpha(x,y)z:=T_\alpha(x,y,z) - T_\alpha(y,x,z) = 
T^+(x,\alpha y,z) - T^+(y,\alpha x,z)
\end{equation}
defines a {\em Lie triple system}.
Thus, given some Jordan pair $(V^\pm,T^\pm)$, its structure variety parametrizes Lie triple
structures on $V^+$. Such families of Lie triple systems are called {\em homotopes}
and have been investigated in Part I of this work (\cite{BeBi}).
 When $\alpha$ is an invertible endomorphism, we speak of {\em isotopes};
in the examples from Part I, the isotopes belong to semisimple or reductive symmetric spaces,
and
when $\alpha$ becomes singular, these spaces contract to degenerate
homotopes, until, for the most singular value $\alpha = 0$, we get the flat space $V^+$. 
In this work, we prove two kinds of results:
\begin{enumerate}
\item 
 A {\em classification theorem} implying that the lists given in Part I are indeed complete
 parametrizations of structure varieties of classical symmetric spaces,
 %(since we are interested
 %in contractions, resp.\ deformations, we really need the complete parametrization of the
 %structure variety, and not only a classification up to isomorphism; the latter can be deduced
 %easily from the former, if needed),
 \item 
 results on the {\em structure} of the homotopes:
 we show that all regular  homotopes have a {\em fibered 
structure} -- they are symmetric spaces fibered over non-degenerate symmetric spaces,
having flat fibers, and such that the fibration splits.
\end{enumerate}
\nin
Both items are closely related to each other: the key ingredient in both cases is
the link with {\em inner ideals}, and with the notion of {\em complementation of inner ideals}
introduced by Loos and Neher (\cite{LoN94}).
More precisely, the image $I$ of the endomorphism $\alpha:V^+ \to V^-$ 
is always an {\em inner ideal in $V^-$},
and the kernel $\ker \alpha$ is a flat ideal in $(V^+,R_\alpha)$ viewed as a Lie triple system.
Thus we get a short exact sequence of Lie triple systems
\begin{equation}
\begin{matrix}
0 & \to & \ker \alpha & \to & (V^+,R_\alpha) & \to & I & \to &  0 \end{matrix} .
\end{equation}
On the level of symmetric spaces, this corresponds to a fibration over a base 
belonging to $I$. 
In general, $I$ does not carry a Jordan triple structure, nor does this sequence split.
We introduce a natural notion of {\em regularity of homotopies} and prove that both properties
hold for regular $\alpha$.
For simple real finite-dimensional
Jordan pairs, we show that every $\alpha$ is indeed regular. 

\ssk
The key observation for proving these results is that $\alpha$ is regular if and only if
the inner ideal $I = \im \alpha$ is {\em complemented}, i.e., there is an inner ideal $J$ in $V^+$
such that
$V^+ = J \oplus \Ker(I)$ and $V^- = I \oplus \Ker(J)$.
By results of Loos and Neher (\cite{LoN94}) it is known that, e.g.\ for simple finite-dimensional
real Jordan pairs, every inner ideal admits a complement.
 This opens the way for classification:
one needs, on the one hand, the (known) classification of inner ideals $I$ in $V^-$, and on the other
hand, for given $I$, in order to describe all $\alpha$ with $I=\im \alpha$, we have
to use the classification of  isotopes of $I$, i.e., of the
invertible elements in $\Svar(I)$ (which is also known).
Putting these two together, one can finally show that, for classical Jordan pairs, the
lists given in Part I (\cite{BeBi}, Theorem 4.2) lead to a complete parametrization of 
structure varieties (with one notable exception:
spaces of skew-symmetric matrices $\Asym(n,\K)$ contain a ``degenerate'' class of inner
ideals, called  ``point spaces'' by McCrimmon \cite{McC}, and homotopes corresponding to such
inner ideals are not given by the construction from Part I). 

\ssk
The organization of this paper is as follows: Chapter 1 contains definitions and basic facts
about Jordan pairs and their structure varieties; in Chapter 2 we introduce the notion of
regularity and explain the link with inner ideals; Chapter 3 contains classifications and
some more results on special cases. While for classification we make quite restrictive
assumptions (base field $\K=\R$, finite-dimensionality and simplicity; this framework is most
important for the theory of symmetric spaces), most of the results
from Chapters 1 and 2 are valid very generally (arbitrary Jordan pairs over general base rings).

\ssk
{\bf Notation.} Throughout, $\K$ denotes a commutative base ring in which $2$ and $3$ are  invertible.

\section{The structure variety of a Jordan pair}

A {\em $\Gamma$-graded Lie algebra} is a Lie algebra of the form
$\g = \bigoplus_{\gamma \in \Gamma} \g_\gamma$ with
$[\g_\gamma,\g_\delta]\subset \g_{\gamma + \delta}$ where $(\Gamma,+)$ is an
abelian group. 

\subsection{$\Z/(2)$-graded Lie algebras and Lie triple systems} \label{sec:Liealg}

Recall that $\Z/(2)$-graded Lie algebras $\g= \h \oplus \q$ (i.e.\ $[\h,\h] \subset \h$,
$[\q,\q] \subset \h$, $[\h,\q] \subset \q$)
essentially correspond to
{\em Lie triple systems (LTS)} by considering $\q$ with the triple Lie bracket 
$[X,Y,Z]:=[[X,Y],Z] =:R(X,Y)Z$. Axiomatically, this triple Lie bracket can be characterized by the three properties

\begin{description}
\item[(LT1)] $[X,Y,Z] = - [Y,X,Z]$
\item[(LT2)] $[X,Y,Z]+[Y,Z,X] + [Z,X,Y]=0$
\item[(LT3)] the endomorphism $D:=R(U,V)$ is a derivation of the trilinear product
$[X,Y,Z]$. 
\end{description}

\nin
Every LTS $\q$ is obtained from a $\Z/(2)$-graded Lie algebra: 
we may take for $\g$ the {\em standard imbedding}
$\q \oplus [\q,\q] \subset \q \oplus \Der(\q)$ (see \cite{Lo69}).
Then $\sigma := \id_\h \oplus (-\id_\q)$ is an automorphism of order $2$, and
the pair $(\g,\sigma)$ is called a {\em symmetric pair}. 
See Section 1 of Part I for further remarks on Lie triple systems.

\subsection{Jordan pairs, $3$-graded Lie algebras and polarized Lie triple systems}

A {\em $3$-graded Lie algebra} is a $\Z$-graded Lie algebra 
with $\g_n = 0$ for $\vert n \vert > 1$:
$$
\g = \g_{-1} \oplus \g_0 \oplus \g_1 .
$$
Note that then the odd part $\q := \g_1 \oplus \g_{-1}$ is a LTS, which, moreover, is
{\em polarized}:  the endomorphism
$D:=\id_{\g_1} \oplus (-\id_{\g_{-1}})$ is a polarization of the LTS $\q$, that is,
$D^2 = 1$, and $D$ is a derivation of the LTS commuting with $\h$ in the sense that
$R(X,Y)\, D=D \, R(X,Y)$.  
Conversely, the standard imbedding of a polarized LTS is a
$3$-graded Lie algebra. % (see, e.g., \cite{Be00}, Chapter III). 
If $\g$ is a $3$-graded Lie algebra, the pair
$(V^+,V^-):=(\g_1,\g_{-1})$ together with the trilinear maps
$$
T^\pm : V^\pm \times V^\mp \times V^\pm \to V^\pm, \quad
(u,v,w)\mapsto T^\pm(u,v,w):=[[u,v],w]
$$
is called the {\em associated (linear) Jordan pair}. It satisfies the
identities

\begin{description}
\item[(J1)] $T^\pm(u,v,w) = T^\pm( w,v,u)$
\item[(J2)] 
$T^\pm(u,v,T^\pm(x,y,z)) = T(T(u,v,x),y,z) - T(x,T^\mp(v,u,y),z) + T(x,y,T(u,v,z))$
\end{description}

\begin{dfn}
A {\em (linear) Jordan pair} is a pair of vector spaces $(V^+,V^-)$
with trilinear maps $T^\pm$ satisfying (J1) and (J2).
{\em Homomorphisms of Jordan pairs} are pairs of linear maps
$(f^\pm:V^\pm \to W^\pm)$ such that
$f^\pm (T^\pm(u,v,w))=T^\pm(f^\pm u,f^\mp v,f^\pm w)$.
Obviously,
Jordan pairs form a category.
\end{dfn}

\nin Every Jordan pair is obtained from a 3-graded Lie algebra
in the way just described.

 \begin{dfn}\label{OppositeDef}
For any Jordan pair $(V^+,V^-)$, the {\em opposite pair} is
$(V^-,V^+)$.  An {\em involution} of $(V^+,V^-)$ is an isomorphism
$\phi^\pm:V^\pm \to V^\mp$ onto the opposite pair such that
$\phi^- = (\phi^+)^{-1}$.
\end{dfn}

\subsection{Involutions, Jordan triple systems and the Jordan-Lie functor}

Assume $\g$ is a $3$-graded Lie algebra and
 $\theta$  an automorphism of $\g$ of order two and reversing the
grading. Then $\theta$ restricts to an automorphism of the LTS
$\q = \g_1 \oplus \g_{-1}$ anticommuting with the polarization $D$,
called a {\em para-real form of $(\q,R,D)$},
and the restriction $\g_1 \to \g_{-1}$ defines an involution of the Jordan pair $(V^+,V^-)$.
Conversely, every para-real
form (every involution) induces an automorphism $\theta$ as above.
If $(\g,\theta)$ is as above, the space $V:=\g_1$ with trilinear product
$$
T(x,y,z):= [[x, \theta(y)],z]
$$
is a {\em Jordan triple system (JTS)}, i.e., it satisfies the identities (J1) and (J2) with superscripts
$\pm$ omitted. Every Jordan triple system $(V,T)$ is obtained in this way from an
involutive $3$-graded Lie algebra (see, e.g.,  \cite{Be00}):
Jordan triple systems are equivalent to Jordan pairs with involution, and
 the Jordan pair $(V^+,V^-)$ obtained by forgetting the involution
is then called the {\em underlying Jordan pair of $(V,T)$}.
Note that in general there are many grade-reversing involutions, and
 there is no canonical choice of one of them; therefore it is often conceptually
clearer to work with Jordan pairs and not with triple systems (as we will do below).

\begin{lemma}\label{LieJordanLemma}
If $(V,T)$ is a Jordan triple system, then $R:=R_T$ with
$$
R_T(x,y,z):=T(x,y,z) - T(y,x,z)
$$
defines a LTS on $V$.
This defines a functor $T \mapsto R_T$ from the category of JTS to the one of LTS, 
called the {\em Jordan-Lie functor}.
\end{lemma}

\Pf Direct computation -- see \cite{Mey73} or \cite{Be00}, Lemma III.2.6. One may note also that
the LTS $(V,R_T)$ is isomorphic to the fixed point space under $\theta$ of the para-real form of
$(\q,R,D)$ mentioned above.
\EPf

\subsection{Definition of the structure variety}\label{sec:Svar}

Let $V=(V^+,V^-)$  be a Jordan pair.  We define the {\em quadratic operators} for $x \in V^\pm$
by
$$
Q^\pm(x):V^\mp \to V^\pm, \quad y \mapsto \frac{1}{2} T^\pm(x,y,x).
$$

\begin{dfn}\label{SvarDef}
The {\em structure variety (of $V^+$)}
is the subset of $\Hom(V^+,V^-)$ defined by Equation (\ref{SvarDefEqn}).
An element $\alpha \in \Svar(V^+)$ is called a {\em homotopy of the Jordan pair}. 
Equivalently (by polarization), a linear map $\alpha$ is a homotopy iff
for all $u \in V^+$,
$$
Q^-(\alpha u)=\alpha \circ Q^+(u) \circ \alpha.
$$
The set of invertible elements in $\Svar(V^+)$ is denoted by $\Svar(V^+)^\times$.
For a JTS $(V,T)$, the structure variety $\Svar(V)$
is defined in the same way, omitting the superscrips $\pm$.
In other words, the structure variety of a JTS is the structure variety of the underlying Jordan pair. 
\end{dfn}

\nin {\bf Remarks.}
(1) If $\alpha$ is an invertible element in $\Svar(V^+)$, then, by change of variables
$v':=\alpha v$, we see that $(\alpha,\alpha^{-1})$ is an involution. Conversely, 
if $(\phi^+,\phi^-)$ is an involution, then $\phi^+ \in \Svar(V^+)^\times$.
In this sense $\Svar(V^+)^\times$ corresponds to the space of involutions, and thus
$\Svar(V^+)$ can be seen as a sort of completion of this space.

\ssk \nin (2)
 The {\em fundamental formula} from Jordan theory
$Q(Q(x)y)=Q(x)Q(y)Q(x)$ says that 
all  quadratic operators 
$Q^-(x):V^+ \to V^-$,
where $x \in V^-$, belong to $\Svar(V^+)$. 

\ssk \nin
(3) Obviously, $\Svar(V^+)$ is stable under multiplication by scalars and contains $0$.

\ssk \nin
(4) There is a similar
definition of $\Svar(V^-)$. From a categorial point of view, one may
interpret the structure variety
as the {\em space of self-adjoint structural transformations from
a Jordan pair to its opposite pair}. Here are the relevant definitions:

\begin{dfn}\label{StructuralTrfDef}
A {\em structural transformation} between two Jordan pairs
$(V^+,V^-)$, $(W^+,W^-)$ is a pair of linear maps 
$(g:V^+ \to W^+, h:W^- \to V^-)$ such that, for all
$u,w \in V^+$ and $ v,z \in W^-$:
$$ 
g T^+_V(u,hv,w)=T^+_W(gu,v,gw), \quad
h T^-_W(v,gw,z)=T^-_V(hv,w,hz) .
$$
Jordan pairs with structural transformations form a category in which composition
of morphisms is defined by $(g,h)(g',h'):=(g \circ g', h' \circ h)$.
In particular,
the structural transformations from $(V^+,V^-)$ to itself form 
a semigroup, called the {\em structure monoid}. By definition,
the {\em structure group of $V$} is the group of invertible elements
of the structure monoid, that is, the group of automorphisms of $(V^+,V^-)$.
A structural transformation $(g,h)$ from $(V^+,V^-)$ to its opposite pair
$(V^-,V^+)$ is called {\em self-adjoint} if $g=h$. 
\end{dfn}
With these definitions, $\alpha$ belongs to $\Svar(V^+)$ if and only if $(\alpha,\alpha)$ 
is a structural transformation to the opposite pair. Note also that, if
$(f^+,f^-)$ is an isomorphism in the ``usual'' sense, then
$(f^+,(f^-)^{-1})$ is an isomorphism in the new (``structural'') category.

\begin{lemma}\label{StructureLemma}
If $\alpha \in \Svar(V^+)$, then 
$
T_\alpha(x,y,z):=T^+(x,\alpha y,z)
$
defines a Jordan triple system on $V^+$,
called the {\em $\alpha$-homotope of $(V,T)$}.
\end{lemma}

\Pf
Straightforward computation -- see  \cite{Mey73} or \cite{Be00}, Lemma III.4.5.
\EPf

\nin Applying the Jordan-Lie functor to $T_\alpha$, we get by Equation
(\ref{RHomotopeEqn}) a LTS on $V^+$, also called {\em $\alpha$-homotope of $R$}.
Thus we have a map $\alpha \to R_\alpha$ from the structure variety to the variety of
Lie triple products on $V^+$. 
%As is easily seen, this map is equivariant under the action of the structure group. 

\begin{lemma}\label{TransferLemma}
Assume $\alpha \in \Svar(V^+)$ and $(g,h):(W^+,W^-) \to (V^+,V^-)$ a structural transformation.
Then $h \circ \alpha \circ g:W^+ \to W^-$ belongs to $\Svar(W^+)$, and
$$
g:(W^+,T_{h \alpha g}) \to (V^+,T_\alpha), \quad v \mapsto gv
$$
is a JTS-homomorphism (and hence also a LTS-homomorphism $R_{h\alpha g} \to R_\alpha$).
\end{lemma}

\Pf
By definition of composition of morphisms,
$(h,g) (\alpha,\alpha) (g,h)=(h\alpha g,h \alpha g)$, and this is a structural transformation of
$W$. Hence $h \alpha g:W^+ \to W^-$ is self-adjoint structural, i.e., it belongs to $\Svar(W^+)$.
Since
$g T_{h \alpha g}(u,v,w)=
g T(u,h \alpha g v,w) =
T(gu,\alpha gv,gw) =
T_\alpha(gu,gv,gw)$, it follows that $g$ induces a JTS-homomorphism as claimed.
\EPf

\begin{cor}
The structure group of the Jordan pair $V$ acts on the structure variety, and 
if $\alpha$ and $\alpha'$ belong to the same orbit under this action, then the JTS
$T_\alpha$ and $T_{\alpha'}$ are isomorphic.
\end{cor}

\Pf
This is the special case $V=W$ and $(g,h)$ invertible of the lemma.
\EPf

For a classification of homotopes up to isomorphy it thus suffices to
consider structure group orbits in $\Svar(V^+)$. 
Note that the structure group contains all pairs $(r \id, r \id)$ with $r \in \K^\times$, and hence
it follows that $\alpha$ and $r^2 \alpha$ for $r \in \K^\times$
 are conjugate under
the structure group. In general, if $\K=\R$,
  $\alpha$ and $-\alpha$ will not
be conjugate to each other (they are {\em $c$-duals} of each other, cf.\ remarks in
Part I). 

\begin{cor}
Assume $\alpha \in \Svar(V^+)$ and
$\beta \in \Svar(V^-)$. Then $\beta \alpha \beta$ belongs to $\Svar(V^+)$,
and $\alpha \beta \alpha$ belongs to $\Svar(V^-)$. 
\end{cor}

\Pf
This is the special case  $(W^+,W^-)=(V^-,V^+)$ and $g=h=\beta$
(resp.\ $g=h=\alpha$) of the lemma.
\EPf

\section{Fibration: Inner ideals, kernels and complementation}

The symmetric space $M_\alpha$ corresponding to the LTS $R_\alpha$
has a {\em fibered structure}: the base corresponds to the image of $\alpha$ and the
fiber to the kernel of $\alpha$. In this chapter we investigate these important features
and show that they are also the key for proving classification results. It turns out that 
classification is made possible by the fact that, in suitably ``regular'' situations, the above
mentioned fibering splits, and that such splittings correspond, in a Jordan theoretic
language, to {\em complementation of inner ideals} in the sense of \cite{LoN94}.

\subsection{Inner ideals and their kernels}

We are going to describe some properties of the image
$\im \alpha \subset V^-$ and the kernel $\ker \alpha \subset V^+$ of a
homotopy $\alpha :V^+ \to V^-$.

\begin{lemma} Assume $\alpha \in \Svar(V^+)$. Then
$\ker \alpha$ is an ideal in the Jordan triple system $(V^+,T_\alpha)$, and hence it is also
an ideal in the Lie triple system $(V^+,R_\alpha)$.
\end{lemma}

\Pf
Let $K:=\ker \alpha$; in order to check that
$T_\alpha(K,V^+,V^+) + T_\alpha(V^+,K,V^+)+T_\alpha(V^+,V^+,K) \subset K$,
it suffices to check that
$\alpha T_\alpha(u,v,w)=0$ whenever one of $u,v,w$ belongs to $K$.
But this follows immediately from the relation $\alpha T_\alpha(u,v,w) =
\alpha T(u,\alpha v,w) = T (\alpha u,v,\alpha w)$.
\EPf

The image $\im \alpha$ is in general not even a ``sub-JTS'' of $V^-$
(so far there is no JTS-structure on $V^-$). 

\begin{dfn} Let $(V^+,V^-)$ be a Jordan pair.
 An {\em inner ideal (in $V^-$)}  is a subspace
$I \subset V^-$ such that 
$$
T^-(I,V^+,I) \subset I.
$$
\end{dfn}

\begin{ex}
Images of structural transformations are inner ideals: assume
$(g,h):(V^+,V^-)\to (W^+,W^-)$ is a structural transformation, then
the  $\im h$ is an inner ideal in $W^-$ since
$$
T^-(hV^-,W^+,hV^-) \subset h T^-(V^-,gW^+,V^-) \subset \im h.
$$
In particular, for $\alpha \in \Svar(V^+)$,
the image $\im \alpha$ is an inner ideal in $V^-$. 
\end{ex}

\begin{dfn} \label{KernDef}
Let $I$ be an inner ideal in $V^-$.  We define, following \cite{LoN94},
the {\em kernel of $I$} to be the subspace of $V^+$ given by
$$
\Ker I := \{ x \in V^+ \mid \, T^-(I,x,I)=0 \} = \bigcap_{y \in I} \ker Q^-(y) .
$$
Note that the kernel of an inner ideal and the inner ideal itself live in different spaces
($V^+$, resp.\ $V^-$).
We recall also that a Jordan pair
$(V^+,V^-)$ is called {\em non-degenerate}
if $T^\mp(x,V^\pm,x)=0$ implies $x=0$.
\end{dfn}

\begin{ex}\label{KernelExample}
The kernel of an image $\im g$ of a structural transformation $(g,h)$ is the kernel of 
$h$:
$$
\Ker (\im g) = \ker h .
$$
 More precisely,
if $(g,h):(V^+,V^-) \to (W^+,W^-)$ is structural, let $I:=\im g$ and
$K:=\ker h$; then
$T^-(gu,k,gv)=g T^-(u,hk,v)=0$ for all $k \in K$, whence
$\ker h \subset \Ker(\im g)$. 
The converse holds if $V$ is non-degenerate
(\cite{LoN94}, Lemma 1.6 (b)).
\end{ex}

\subsection{Regularity and pseudo-inverses}\label{REG}

An inner ideal $I$ does in general not carry a Jordan triple system structure
(we can only say that $(V^+,I)$ is a Jordan pair), and therefore it does not make
sense to say that $T_\alpha /\ker \alpha$ and $\im \alpha$ be isomorphic as JTS
(and thus as LTS). However, this will be the case if $\alpha$ is {\em regular} in the
sense to be explained now.

\begin{dfn}\label{RegDef}
Let $(V^+,V^-)$ be a Jordan pair and $\alpha \in \Svar(V^+)$ (hence
$\alpha:V^+ \to V^-$).
We say that $\alpha$ is {\em regular} if there is an element
$\beta \in \Svar(V^-)$ (so $\beta:V^- \to V^+$) such that
$$
\alpha \circ \beta \circ \alpha = \alpha \, .
$$
(This implies that $\alpha$ is regular in the
Jordan pair $(\Hom(W^+,W^-),\Hom(W^-,W^+)$, see \cite{Lo75}.)
The element $\beta$ will then be called a {\em pseudo-inverse of $\alpha$}.
\end{dfn}

\nin In the situation of the definition, Lemma \ref{TransferLemma} implies that
$$
\alpha : (V^+,T_{\alpha \beta \alpha}) = (V^+,T_\alpha) \to (V^-,T_\beta)
$$
is a homomorphism of Jordan triple systems (and hence also a homomorphism of
LTS from $R_\alpha$ to $R_\beta$). 

%\Pf
%$\alpha T_\alpha(u,v,w) = \alpha T(u,\alpha v,w) = \alpha T(u,\alpha \beta \alpha v,w) 
%= T(\alpha u, \beta \alpha v, \alpha w) = T_\beta(\alpha u, \alpha v, \alpha w) $
%\EPf

\begin{ex} Assume  $(V,T)$ is a JTS, hence $V=V^+=V^-$. 
If $\alpha$ is a tripotent  ($\alpha^3 = \alpha$) element of $\Svar(V)$, then it is regular
(take $\beta = \alpha$), and hence
$\alpha:T_\alpha \to T_\alpha$ is a JTS-homomorphism.
This holds  in particular if $\alpha^2 = \alpha$ or $(-\alpha)^2 = -\alpha$.
However, not every regular element is tripotent in this sense.
%
%for instance, if $\alpha' = h \circ \alpha \circ g$ for $(g,h)$ in the structure group,
%the JTS $T_\alpha$ and $T_{\alpha'}$ will be isomorphic, hence are both regular,
% but $\alpha'$ needs no longer be tripotent. 
\end{ex}

If the set of pseudo-inverses of $\alpha$ is not empty, then it is the intersection of $\Svar(V^-)$
with an affine space: the
difference $\gamma=\beta - \beta'$ of two pseudo-inverses satisfies the linear equation
$\alpha \gamma \alpha = 0$.  There is no canonical choice of pseudo-inverse, but once
we have fixed a choice, we may ``improve'' it in a canonical way: 
by straightforward computation one proves the following

\begin{lemma}
If $\beta$ is a pseudo-inverse of $\alpha$, then the endomorphism
$\beta':=\beta \alpha \beta : V^- \to V^+$ belongs again to $\Svar(V^-)$, is again
a pseudo-inverse of $\alpha$, and it satisfies moreover
$$
\beta' \circ \alpha \circ \beta' = \beta' \, .
$$ 
\end{lemma} 

\begin{dfn}
Let $(W^+,W^-)$ be a pair of $\K$-modules and $\alpha :W^+ \to W^-$, $\beta:W^- \to W^+$
linear maps.
We say that a pair $(\alpha,\beta) \in \Hom(W^+,W^-) \times \Hom(W^-,W^+)$ is {\em idempotent}
if $\beta \circ \alpha \circ \beta = \beta$ and $\alpha \circ \beta \circ \alpha = \alpha$.
(Equivalently, $(\alpha,\beta)$ is an idempotent in the
Jordan pair $(\Hom(W^+,W^-),\Hom(W^-,W^+)$, see \cite{Lo75}.)
\end{dfn}

\begin{lemma} %see http://fr.wikipedia.org/wiki/Pseudo-inverse
\label{DecompositionLemma}
Assume 
 the pair $(\alpha,\beta)$ is idempotent. Then we have direct sum decompositions
\begin{equation} \label{direct!}
W^+ = \im \beta \oplus \ker \alpha, \quad
W^- = \im \alpha \oplus \ker \beta,
\end{equation}
and $\beta\vert_{\im \alpha}:\im \alpha \to \im \beta$ is a linear isomorphism with inverse
$\alpha\vert_{\im \beta}:\im \beta \to \im \alpha$.
Conversely, if a pair $(\alpha,\beta)$ has these properties, then it is idempotent.
\end{lemma}

\Pf
Direct check (cf.\  \cite{BeL08}, Lemma A.8; see also http://fr.wikipedia.org/wiki/Pseudo-inverse).
\EPf

\begin{thm} \label{Complement-Theorem}
If $\alpha \in \Svar(V^+)$ is regular, then $\im \alpha$ carries a JTS structure such that
$$
0 \to \ker \alpha \to (V^+,T_\alpha) \to \im \alpha \to 0
$$
is a split exact sequence of JTS. A similar statement holds for the corresponding LTS.
\end{thm}

\Pf
Choose $\beta$ such that $(\alpha,\beta)$ is idempotent.
According to the preceding lemma, $\alpha$ and $\beta$ induce mutually inverse
isomorphisms between $I=\im \alpha$ and $J:=\im \beta$. 
Therefore $I$ becomes a JTS with triple product
$(u,v,w) \mapsto T^-(u,\beta v,w)$, and
$\alpha:(V^+ ,T_\alpha) \to I$ becomes a homomorphism which splits via $\beta:I \to V^+$.
\EPf

By general functoriality of the construction of symmetric spaces associated to Jordan structures
(\cite{Be02}), we deduce that the exact sequence from the theorem lifts to the space
level: $\alpha$ induces a homomorphism
from $M_\alpha$ to the symmetric space $B_\alpha$ belonging to the base
$\im \alpha$; we have a fibration of symmetric spaces
$$
F_\alpha \to M_\alpha \to B_\alpha,
$$
where the fiber $F_\alpha$ is a flat symmetric space, and there exists a splitting 
$B_\alpha \to M_\alpha$. 
The LTS $(V^+,R_\alpha)$ of $M_\alpha$ splits as
$V^+ = K \oplus J$ with abelian ideal $K$ and $J:=\im \beta$. Moreover, we have

\begin{lemma}
In the LTS $(V^+,R_\alpha)$ with ideal $K=\ker \alpha$ we have
$
R_\alpha(K,K)V^+ = 0 .
$
\end{lemma}

\Pf
$T_\alpha(u,v,w) - T_\alpha(v,u,w) =
T(u,\alpha v,w) - T(v,\alpha u,w) = 0$
whenever $u,v \in K = \ker \alpha$.
\EPf

However, $R_\alpha(K,J)K$ is in general not zero, and hence the fibration $M_\alpha \to B_\alpha$
is in general not a {\em symmetric bundle} in the sense of \cite{BeD09}.
%
%\bigskip 
%[Some related
%Questions and remarks, to be removed or explained: is $R(K,I)$ solvable, $[R(K,I),R(K,I)]$ nilpotent? t
%note $R_\alpha(K,J)V^+ \subset K$ since it is an ideal 
%
%but $R_\alpha(K,J)J \subset ?$ or $R_\alpha(K,J)K=0$ ? if non-deg ? 

%note $R(F,V)=T(F,\alpha V)- T(V,\alpha F)=T(F,\alpha V)=T(F,B)=R(F,B)$

\subsection{Complementation of inner ideals}

Recall from \cite{LoN94} that an inner ideal $J \subset V^+$ is called a {\em complement}
of an inner ideal $I \subset V^-$ if 
$$
V^- = I \oplus \Ker J, \quad V^+ = J \oplus \Ker I \, ,
$$
and that Jordan pair $(V^+,V^-)$ is
 {\em complemented} if every inner ideal both in $V^+$ and $V^-$ admits a
complement. Complemented Jordan pairs have been characterized in \cite{LoN94};
we will use their result in the proof of the following theorem, but we will give also a more
elementary, independent proof for the positive real case in the next subsection.

\begin{thm}\label{RegularityTh}
Assume $(V^+,V^-)$ is a non-degenerate Jordan pair  that
admits an {\em anisotropic involution $\tau$} (this means that $T^+(x,\tau(x),x)=0$ implies $x=0$;
one says also that the JTS $(V,T)$ with $T(x,y,z)=T^+(x,\tau(y),z)$ is anisotropic).
\begin{description}
\item[(1)] If $\alpha \in \Svar(V^+)$ is regular, then 
the inner ideal $I:=\im \alpha$ is complemented.
A complement is given by   $J:=\im \beta \subset V^+$, where $\beta \in \Svar(V^-)$ is such
that $(\alpha,\beta)$ is idempotent.
\item[(2)] Assume that $(V^+,V^-)$ is simple and finite-dimensional  over a field $\K$.
Then every inner ideal $I \subset V^-$ is complemented, and there exists
regular $\alpha \in \Svar(V^+)$ such that $I = \im \alpha$. Moreover, every element of
the structure variety is regular, and after fixing an anisotropic involution,
 all elements of $\Svar(V^+)$ having same 
image and kernel as $\alpha$ are naturally parametrized by $\Svar(I)^\times$. 
\end{description}
All assumptions hold, in particular, if $(V^+,V^-)$ is a simple, finite-dimensional
Jordan pair over $\K=\R$; we can then choose for $\tau$ its {\em Cartan-involution},
i.e., the involution coming from a Cartan-involution of the corresponding $3$-graded Lie
 algebra reversing the grading. 
\end{thm}

\Pf
(1) From Lemma \ref{DecompositionLemma} we get the decomposition
\begin{equation}\label{compl}
V^+ = \im \beta \oplus \ker \alpha, \quad
V^- = \im \alpha \oplus \ker \beta .
\end{equation}
Since  our assumption implies that $(V^+,V^-)$ is non-degenerate, we have
$\ker \alpha = \Ker(\im \alpha)$ (cf.\ Example \ref{KernelExample}) and
$\ker \beta = \Ker(\im \beta)$. This shows that $\im \beta$ is a complement of
$\im \alpha$. 

\ssk
(2) We fix for the moment
 an anisotropic involution $\tau:V^+ \to V^-$, giving rise to the JTS $(V,T)$.
 By definition, this JTS is then anisotropic, and it is Artinian (i.e., descending chains of
 inner ideals become stationary) since it is finite-dimensional. Thus
\cite{LoN94}, Theorem 6.7 implies that $I$ is complemented in the Jordan triple sense, that is,
$J:=\tau(I) \subset V^+$ is a complement of $I \subset V^-$. 
Let $\pi^\pm : V^\pm \to V^\pm$ be the projections onto $I$, resp.\ $J$ with kernel
$\Ker J$, resp.\ $\Ker I$; then $\tau \circ \pi^+ = \pi^- \circ \tau$.
As shown in \cite{LoN94}, Section 3.5, the pair of projections
 $(\pi^+,\pi^-)$ is a structural transformation, 
and hence also the pair
$$
(\tau,\tau)  (\pi^+,\pi^-) = (\tau \circ \pi^+,\pi^- \circ \tau) =
(\tau \circ \pi^+, \tau \circ \pi^+) =: (\alpha_0 , \alpha_0 )
$$
is structural (mind the definition of composition in the structural category). Similarly, the pair
$(\beta_0,\beta_0)$ with $\beta_0 : = \tau^{-1} \circ \pi^-$ is structural, whence
$\alpha_0 \in \Svar(V^+)$, $\beta_0 \in \Svar(V^-)$. Moreover,
$$
\alpha_0 \beta_0 \alpha_0 = \tau \circ \pi^+ \circ \tau^{-1} \circ \pi^- \circ \tau \circ \pi^+
=
\tau \circ (\pi^+)^3 = \alpha_0
$$
and similarly $\beta_0 \alpha_0 \beta_0 = \beta_0$. By definition, the image of $\alpha_0$ is $I$ and
its kernel is $\Ker I$, proving that every inner ideal $I$ is the image of a regular $\alpha_0$.
Still fixing the anisotropic involution, we have the following lemma describing all 
elements $\alpha \in \Svar(V^+)$ having same image and kernel as $\alpha_0$:

\begin{lemma}\label{SvarILemma}
Let $I$ be a complemented inner ideal in a JTS $V$
and let $\pi:V \to I$ and $\iota:I \to V$ be  projection and injection corresponding to
the decomposition $V = I \oplus \Ker I$.
Then $(\pi,\iota):(V,V) \to (I,I)$ and
$(\iota,\pi):(I,I) \to (V,V)$ are structural transformations, and we obtain two maps
\begin{eqnarray*}
\phi: \Svar(I) \to \Svar(V), & \quad & \gamma \mapsto \alpha := \iota \circ \gamma \circ \pi 
= \begin{pmatrix} \gamma & 0 \cr 0 & 0 \end{pmatrix}
\cr
\psi: \Svar(V) \to \Svar(I), & \quad & \alpha = 
\begin{pmatrix} \gamma & \gamma' \cr \gamma'' & \gamma''' \end{pmatrix} 
\mapsto \gamma = \pi \circ \alpha \circ \iota
\end{eqnarray*}
(matrices taken w.r.t.\ the decomposition $V = I \oplus \Ker I$)
such that $\psi \circ \phi = \id_{\Svar(I)}$. In particular, $\phi$ is injective and $\psi$ is
surjective, and the image of $\phi$ is precisely the set of elements $\alpha \in \Svar(V)$ such
that $\im \alpha \subset I$ and $\Ker I \subset \ker \alpha$. 
In particular, the conditions $\im \alpha = I$ and $\ker \alpha = \Ker I$ hold if and only if
$\alpha = \phi(\gamma)$ with {\em invertible} $\gamma \in \Svar(I)$. 
\end{lemma}

\Pf 
As mentioned above, the pair of projections $(\pi^+,\pi^-)$ is structural.  Therefore
$\pi T(u,\iota v,w) = T(\pi u,v,\pi w)$ for 
$u,w \in V$, $v \in I$, hence $(\pi,\iota):(V,V) \to (I,I)$ is structural.
For $(\iota,\pi)$: let $u,w \in I$, $v = v_I + v_K$ with $v_I \in I$, $v_K \in \Ker (I)$, then
$$
\iota T(u,\pi v,w) = T(u,v_I,w) = T(u,v_I + v_K,w) = T(u,v,w) = T(\iota u,v,\iota w) .
$$
It follows that the pair
$$
(\iota,\pi) (\gamma ,\gamma)(\pi,\iota) = (\iota \circ \gamma \circ \pi, \iota \circ \gamma \circ \pi) =
(\alpha,\alpha) : (V,V) \to (V,V) 
$$
is structural, and hence $\alpha \in \Svar(V)$, so $\phi$ is well-defined.
The proof showing that $\psi$ is well-defined is  similar. 
Using the ``matrix description'', the remaining statements follow immediately.
\EPf

We now finish the proof of part (2) of the theorem:
let $\alpha \in \Svar(V)$ be arbitrary and put $I:=\im \alpha$, $J:=\tau(I)$.
Then $J$ is a complement of $I$, hence $\im \alpha = \im \alpha_0$ and
$\ker \alpha = \ker \alpha_0$ with regular $\alpha_0$. The lemma shows
that $\alpha_0 = \phi(\id_I)$ and $\alpha = \phi(\gamma)$ for some
$\gamma \in \Svar(I)^\times$. Then it is easily checked that $\alpha$ is regular with
pseudo-inverse $\beta :=\phi(\gamma^{-1})$.

\ssk
In order to prove the final assertion of the theorem,
let us show that the Cartan-involution is anisotropic:
indeed, the JTS corresponding to this involution is {\em positive} (cf.\ next subsection).
 Let $x \in V$ and fix a Jordan frame $e_1,\ldots,e_r$ such
that $x = \sum_i r_i e_i$ with $r_i \geq 0$. Then $T(x,x,x)=\sum_i r_i^3$,
hence $T(x,x,x)=0$ implies $r_i = 0$ for all $i$, hence $x=0$.
%
%geometrically: there are no geodesics that are affine lines)
\EPf

The preceding results lead to a method of classification of elements of the structure
variety of a finite-dimensional real simple Jordan pair $(V^+,V^-)$ in two steps:
Part (2) of the theorem says that the map
$$
\Svar(V) \to \cX:= \{ \mbox{inner ideals in } V^- \}, \quad \alpha \mapsto I := \im (\alpha)
$$
is surjective with fiber over $I$ being in one-to-one correspondence with
$\Svar(I)^\times$. Therefore $\Svar(V)$ can be described by combining the classification
of inner ideals $I$ with the one of elements of $\Svar(I)^\times$. 
Before coming to classification, we give the promised elementary proof for the real case.
It relies on the fact that kernels of inner ideals can be expressed in terms of
 {\em orthogonal
complements} with respect to a suitable bilinear form.

\subsection{The trace form}

In this section,
let $(V^+,V^-)$ be a finite-dimensional Jordan pair over a field $\K$. The 
{\em trace form} is the bilinear form defined by
$$
V^+ \times V^- \to \K, \quad (u,v) \mapsto \langle u,v \rangle := \Tr (T^+(u,v,\cdot) )
$$
If $V$ is a Jordan triple system, then the same form defines a bilinear form
on $V$, again called {\em trace form}.
It is known and easy to show that, if this form is non-degenerate, then it is symmetric, and that
it is {\em invariant} in the sense that the adjoint of $T(x,y,\cdot)$ is $T(x,y,\cdot)$:
\begin{equation}
\langle T(x,y,u),v \rangle = \langle u,T(y,x,v) \rangle .
\end{equation}
If $\K=\R$, we say that a JTS $(V,T)$ is {\em positive}, if the trace form is positive definite,
and {\em negative} or {\em compact} if it is negative definite.
Of course, replacing $T$ by $-T$ makes both notions equivalent. 

\begin{thm}
Assume $(V,T)$ is a real, positive JTS, and denote by  $\perp$ the
  orthogonal complement with respect to the trace form.
\begin{description}
\item[(1)]
Assume
 $I$ is an inner ideal in $V$. Then
$$
\Ker I  = I^\perp \, .
$$
In particular,
the inner ideal $I$ is complemented in
the Jordan triple sense: $V = I \oplus \Ker I$. 
\item[(2)] 
Any $\alpha \in \Svar(V)$ is symmetric with respect 
to the trace form.
\end{description}
\end{thm}

\Pf
(1) We will use the following lemma 

\begin{lemma}
Let $(V,T)$ be a positive JTS. Then every sub-JTS $F \subset V$ is again positive, and
satisfies in particular $T(F,F,F)=F$.
\end{lemma}

\nin Before proving the lemma, let us use it to establish the theorem: 
for all $u,v,w \in I$ and $k \in \Ker I$,
$$
\langle T(u,v,w),k \rangle = \langle w,T(v,u,k) \rangle =
\langle w,T(k,u,v) \rangle = \langle T(u,k,w),v \rangle = 0 ,
$$
whence, using the lemma,
$\Ker I \subset T(I,I,I)^\perp = I^\perp$.
Let us prove that $I^\perp \subset \Ker I$, i.e., $T(I,I^\perp,I)=0$.
By definition of inner ideals, we have $T(I,I^\perp,I) \subset T(I,V,I) \subset I$.
On the other hand, using again the lemma, 
$$
0 = \langle I, I^\perp \rangle = \langle T(I,I,I),I^\perp \rangle =
\langle T(I,I^\perp,I),I \rangle ,
$$
showing that $T(I,I^\perp,I) \subset I^\perp$, whence
$T(I,I^\perp,I) \subset (I \cap I^\perp)=0$, and (1) is proved.

(2)
If $\alpha \in \Svar(V)$ is invertible, then, writing $T(u,v):=T(u,v,\cdot)$, 
$$
T(u,\alpha v) = \alpha \circ T(\alpha u,v) \circ \alpha^{-1},
$$
and taking traces we get $\langle u,\alpha v \rangle = \langle \alpha u,v \rangle$.
Now consider $\alpha \in \Svar(V)$, not necessarily invertible.
Using Lemma \ref{SvarILemma}, we write 
$\alpha = \iota \circ \gamma \circ \pi$, where, according to Part (1), 
$\pi:V \to I$ is the orthogonal projection onto $I:=\im (\alpha)$ with kernel $\Ker I = I^\perp$.
Therefore the adjoint operator of $\pi$ is $\iota$.
Using this, and that $\gamma \in \Svar(I)$ is invertible, whence
$\gamma^*=\gamma$ by the preceding remark, it follows that
$$
\alpha^* = (\iota \circ \gamma \circ \pi)^* =
\pi^* \circ \gamma^* \circ \iota^* = \iota \circ \gamma \circ \pi = \alpha ,
$$
proving (2).
It remains to prove the lemma.
Recall that a real JTS is called {\em Euclidean} if it is endowed with an invariant 
Euclidean scalar product: $(T(x,y)u,v) = (u,T(y,x)v)$.
By restriction of the trace form (that is, $(u,v):=\langle u,v \rangle_V$), 
every subsystem $F$ in a positive JTS is Euclidean.
Let us show that $F$ is again positive, i.e.,
$\langle x,x \rangle_F > 0$ for all non-zero $x \in F$.
To this end, using the spectral theorem in $V$, we decompose
$x = \sum_{i=1}^k r_i e_i$ with $r_i > 0$ and orthogonal non-zero idempotents $e_1,\ldots,e_k$
which can be expressed as certain linear combinations of (odd) powers of $x$;
in particular they belong again to $F$ since $F$ is a subsystem. It follows that
$$
\langle x, x \rangle_F = \sum_i r_i^2 \langle e_i,e_i \rangle_F =
\sum_i r_i^2 \langle T(e_i,e_i,e_i),e_i \rangle_F
$$
Now,
a result due to Backes (see Theorem 1 page 268 of \cite{Backes83})---based on an earlier one of Koecher--- asserts that for all $u,v$ in the Euclidean JTS $F$ the inequality
\begin{equation}\label{BACKES}
\Tr\,S(S(u,v)u,v)\;\geq\;0
\end{equation}
holds with $S(u,v)\;:=\;\frac{1}{2}(T(u,v)+T(v,u))\in\End(F)$. Moreover, 
 equality holds in (\ref{BACKES}) if and only if $S(u,v)=0$.
This implies that $\langle x,x \rangle_F \geq 0$, and that
 $\langle x,x \rangle_F = 0$ iff
$r_i^2 S(e_i,e_i)=0$ for all $i$, iff 
$r_i^2 S(e_i,e_i)e_i=0$ for all $i$, iff
$r_i^2 =0$ for all $i$, iff $x=0$. Thus $F$ is positive. 

Let us show that $T(F,F,F)=F$. Decompose the positive JTS $F$ as $F = U \oplus U^\perp$
with $U = T(F,F,F)$. Then $T(U^\perp,U^\perp,U^\perp)=0$
(since $0=( U,U^\perp ) = (T(F,F,F),U^\perp)= (F, T(F,F,U^\perp)) \supset (F,T(U^\perp,U^\perp,
U^\perp)$), hence $U^\perp$ is a positive JTS with zero product, hence it is zero and
$U = F$.
%
%{\bf Here another
%twist of the argument:}
%Assume $(V,T)$ is Euclidean. Then the scalar product is an invariant scalar product on
%the LTS $R_T$, hence the corresponding symmetric space $M$ is Riemannian. More precisely,
%its Hermitian complexification $M_{h\C}$ is a Hermitian symmetric space of non-compact type
%(\cite{Be00}, Prop.\ V.5.2). The de Rham decomposition of $M_{h\C}$ into irreducible components
%yields a decomposition into simple Hermitian factors, corresponding to a decomposition of
%$(V,T)$ into simple JTS. Now anisotropy guarantees that no factor of type $\R$ with trivial
%product appears in this decomposition, and hence $(V,T)$ is a direct product of
%simple positive JTS and hence is itself positive.
\EPf

For semisimple Jordan {\em algebras}, 
statement (2) is essentially already contained in \cite{Ri70}, Theorem 3 (with a very different proof).
In loc.\ cit.\, Corollary 3.1, 
it is also remarked that this result implies that,  if $(g,h)$ is a structural transformation
of a real semisimple Jordan algebra, then
$h=g^*$ (adjoint operator); the same remark applies in the present context.

\begin{cor}
Assume $(V^+,V^-)$ is a finite-dimensional semisimple Jordan pair over $\R$. Then:
\begin{description}
\item[(1)]
For any an inner ideal $I$ in $V^-$, its kernel $\Ker I \subset V^+$ is the orthogonal complement
$I^\perp$ of $I$ with respect to the trace form, and $I$ admits a complement. 
\item[(2)] 
Any $\alpha \in \Svar(V^+)$ is symmetric with respect 
to the trace form.
\end{description}
\end{cor}

\Pf
It is known that every semi-simple real finite-dimensional Jordan pair admits a
positive involution $\theta$ (the Cartan-involution, see Theorem \ref{RegularityTh}),
 that is, the JTS $V=V^+$
with $T(x,y,z):=T^+(x,\theta y,z)$ is positive. Applying Part (1) of the theorem to this positive
JTS gives Part (1) of the corollary: indeed, both $\Ker I$  and $I^\perp$ are defined 
independently of the involution; we have seen that under identification of $V^+$ and $V^-$ under
$\theta$ the corresponding spaces are equal, therefore $\Ker I = I^\perp$ (and this
equality remains true under identification by any involution, be it positive or not). 
Exactly the same argument shows symmetry of $\alpha$ w.r.t.\ the trace form.
\EPf

The properties from the corollary concern the trace form of Jordan {\em pairs},
whereas  positivity (and anisotropy) are properties of an involution, hence of 
triple systems. Positive (more generally, anisotropic)  forms
permit to choose ``simulatenously'' complements for {\em all} inner ideals;
for base fields different from $\R$, such a choice is in general
not possible, even if a Jordan pair is complemented.

\section{Classification}

It goes without saying that in general a classification of the objects introduced in the preceding
chapters is out of reach. The assumptions that we will make are fairly restrictive and are mainly
imposed by geometric applications: we will look at {\em simple finite-dimensional Jordan pairs
over $\K=\R$}. (From a purely algebraic point of view, it would certainly be interesting to classify objects
over more general base fields and by relaxing finite-dimensionality to chain conditions.)
Moreover, in order to keep this work in reasonable bounds, we restrict attention to {\em
classical} Jordan pairs.
In this chapter we present the (known) classifications of such Jordan pairs and triple systems
 and of their inner ideals, which then lead to a complete description  of their
 structure varieties. 

\subsection{Classification of Jordan pairs}\label{JPclass}

{\bf 
I. Simple complex Jordan pairs}
(cf. \cite{Be00}, Table XII.2.1)

\msk
\noindent
\begin{tabular}{llll}
label & $(V^+,V^-)$ &  $T^\pm$ &
$3$-graded Lie algebra \cr
\hline
1. & $(M(p,q;\C),M(q,p;\C)$ & $T^\pm(X,Y,Z)= X Y Z + Z Y X $ & $\sL(p+q;\C)$
\cr
2. & $(\Sym(n;\C),\Sym(n;\C))$ & $T^\pm(X,Y,Z)= X Y Z + Z Y X $ & $\sP(n,\C)$
\cr
3. & $(\Asym(n;\C),\Asym(n;\C))$ & $T^\pm(X,Y,Z)= X Y Z + Z Y X $
 & $\oo(2n,\C)$
\cr
4. & $(\C^n,\C^n)$ & $T^\pm(x,y,z)= \beta(y,z) x + \beta(y,x) z
- \beta(x,z) y$
 & $\oo(\hat \beta)$
\cr
5. & $(\Herm(3,\OOO_\C),\Herm(3,\OOO_\C))$ & $T^\pm(X,Y,Z)= X Y Z + Z Y X $ &
$\e_7$
\cr
6. & $(M(1,2,\OOO_\C),M(1,2,\OOO_\C))$ & 
$T^\pm(X,Y,Z)= X Y Z + Z Y X $ &
$\e_6$
\end{tabular}

\msk \nin
{\bf Remarks.}
(1) In type 4, $\C^n$ is equipped with a non-degenerate and
symmetric bilinear form $\beta$, and $\hat \beta$ is the
form on $\C^{n+2}$ given by
$\hat \beta((z_0,z,z_{n+1}),(w_0,w,w_{n+1}))=
\beta(z,w) + z_0 w_{n+1} + w_0z_{n+1}$.
The product $T(x,y,z)$
may be realized  in the Clifford algebra $Cl(\C^n,\beta)$ 
as follows:
recall  the Clifford relation $x\bullet y:=\frac{1}{2}(xy+yx)= \beta(x,y)1$
and note that $xyz + zyx =  x \bullet (y \bullet z) - y \bullet (x \bullet z) +
(x \bullet y) \bullet z$.
Using this, one gets for $x,y,z$ belonging to the subspace $\C^n$ of
$Cl(\C^n,\beta)$:
$$
xyz + zyx =  x \bullet (y \bullet z) - y \bullet (x \bullet z) +
(x \bullet y) \bullet z =
\beta(y,z) x + \beta(y,x) z
- \beta(x,z) y.
$$

\nin
(2) Types 4, 5 and 6 are included here for convenience, but 
we we will mainly consider the matrix cases 1, 2 and 3 in the sequel.

\ssk \nin
(3) The following lists will be simpler than the ones
given in \cite{Be00} in the sense that we will not have to 
distinguish for Type 1 the cases $p \not=q$ and $p=q$,  nor
for Type 3 the cases $n$ even and $n$ odd (which had to be
done in \cite{Be00} since only {\em invertible} $\alpha$ were considered there). 

\msk \nin
{\bf II. Simple real Jordan pairs}
(cf. \cite{Be00}, Table XII.2.2).
There are two kinds: the complex pairs from above, considered as real, and the
following:

\msk
\noindent
\begin{tabular}{llll}
label & $(V^+,V^-)$ &  $T^\pm$ &
$3$-graded Lie algebra \cr
\hline
1.1 & $(\Herm(n,\C),\Herm(n,\C))$ 
& $T^\pm(X,Y,Z)= X Y Z + Z Y X $ & $\su(n,n)$
\cr
1.2 & $(M(p,q;\R),M(q,p;\R))$ & $T^\pm(X,Y,Z)= X Y Z + Z Y X $ & $\sL(p+q;\R)$
\cr
1.3 & $(M(p,q;\HHH),M(q,p;\HHH))$ & $T^\pm(X,Y,Z)= X Y Z + Z Y X $ &
 $\sL(p+q;\HHH)$
\cr
2.1 & $(\Sym(n;\R),\Sym(n;\R))$ & $T^\pm(X,Y,Z)= X Y Z + Z Y X $ & $\sP(n,\R)$
\cr
2.2 & $(\Herm(n;\widetilde \HHH),\Herm(n;\widetilde \HHH))$ & 
$T^\pm(X,Y,Z)= X Y Z + Z Y X $ & $\so^*(2n)$
\cr
3.1 & $(\Herm(n;\HHH),\Herm(n;\HHH))$ & 
$T^\pm(X,Y,Z)= X Y Z + Z Y X $ & $\sP(n,n)$
\cr
3.2 & $(\Asym(n;\R),\Asym(n;\R))$ & $T^\pm(X,Y,Z)= X Y Z + Z Y X $
 & $\oo(n,n)$
\cr
4. & $(\R^{p,q},\R^{p,q})$ ($p+q=n$) & $T^\pm(x,y,z)=  \beta(y,z) x + \beta(y,x) z
- \beta(x,z) y$ & $\oo(\hat \beta)$
\end{tabular}

\msk \nin
{\bf Remarks.}
(1) The pair with label $i.j$ is a real form of the complex pair with
label $i.$ from the preceding table (with certain restriction of
parameters: for instance, $1.1$ is a real form in the square case,
i.e., $p=q=n$, and quaternionic matrices are real forms for even
parameters). 
For type 4., $(p,q)$ denotes the signature of the form $\beta$.
Note  that $\R^{p,q}$ and $\R^{p',q'}$ are not isomorphic if
$p$ is different from $p'$ and $q'$.

\ssk \nin
(2) Exceptional Jordan pairs and triple systems have been classified by E.\ Neher.

\subsection{Classification of inner ideals}

For the classical Jordan pairs listed above, we describe their inner ideals.
This classification is known: for Jordan algebras see
\cite{McC}; for Jordan pairs see \cite{Ne96} or \cite{DFGG}; for  an elementary proof in case
of rectangular or symmetric matrices,  see \cite{BeL08}, Appendix A. 

\msk
\nin {\bf 1. Rectangular matrices.} 
 $(V^+,V^-)=(M(p,q;\K),M(q,p;\K)) = (\Hom(E,F),\Hom(F,E))$ over a field $\K$:
all inner ideals in $V^-$ are of the form 
$$
I=I_{a,b}:=\{ f: F \to E \mid \, a \subset \ker f, \im f \subset b \}
$$ 
for linear subspaces $a \subset F$,
$b \subset E$. Thus $I$ can be identified with a matrix space
$\Hom(a',b)$ where $a'$ is some vector space complement of $a$ in $F$. 
This ideal is {\em principal} if and only if it is isomorphic to a space of {\em square} matrices, i.e.,
if $\dim b = \dim a'$. 
In the general case, elementary linear algebra shows that the Kernel of $I_{a,b}$ is
$$
\Ker I_{a,b} = \{ h: E \to F \mid \,  h(b) \subset a \} ,
$$
and $J:=I_{b',a'} \subset \Hom(E,F)$ is a complement of $I_{a,b}$ if $a',b'$ are vector space
complements of $a$, resp.\ of $b$. 

\msk \nin {\bf 2. Symmetric and Hermitian matrices.}
 $(V^+,V^-)=(V,V)$, $V = \Sym(n,\K)$ or $\Herm(n,\bF)$. 
Inner ideals are constructed as above, by taking $E=F=\K^n$ or $\bF^n$ and for $f$ self-adjoint maps 
(w.r.t.\ standard bilinear or Hermitian forms). Then the condition $\im f \subset b$ is equivalent
to $b^\perp \subset \ker f$, i.e., inner ideals are given as above with $a = b^\perp$:
$$
I=I_{b}:=\{ f: E \to E \mid \,  f = f^*,   \im f \subset b \} .
$$ 
The inner ideal $I_b$ is alwyas principal, and it can be identified with the matrix space
$\Herm(b)$. The kernel $\Ker I_b$ is characterized by the property $h(b) \subset b^\perp$,
and under the identification $V^+ = V^-$, a complement of $I$ is given by $I$ itself. 

\msk \nin {\bf 3. Skew-symmetric matrices.}  
 $(V^+,V^-)=(V,V)$, $V = \Asym(n,\K)$.
 There are two different kinds of inner ideals:
The first kind is constructed exactly as in the preceding case:
$$
I=I_{b}:=\{ f: E \to E \mid \,  f = - f^*,   \im f \subset b \} 
$$ 
where $b^\perp \subset \ker  f$ follows automatically. 
Kernel and complement are described as above. 
The second type (called {\em point space} in \cite{McC} and \cite{DFGG}) is fairly
special: fix a non-zero vector $u \in E=\K^n$ and let
$$
K_u := \{ f: E \to E \mid \, f=-f^*, \, f(u^\perp) \subset \K u \} 
$$
(where $\perp$ refers to the standard scalar product on $\K^n$). 
In matrix realization, taking $u$ as the last base vector $e_n$ and $u^\perp = 
\mbox{vect}( e_1,\ldots,e_{n-1})$, these are skew matrices having zero upper left
$(n-1)\times(n-1)$-block.
We check that $K_u$ is an inner ideal: the key observation is that, for any skew-symmetric
$g:E \to E$, we have $\langle  gu,u \rangle = - \langle u,gu \rangle = - \langle gu,u \rangle$, whence
$g(\K u) \subset u^\perp$; therefore, for any $f \in K_u$,
$$
(fgf)(u^\perp) \subset fg (\K u) \subset f(u^\perp) \subset \K u,
$$
whence $fgf \in K_u$. 
Observe also that $K_u = \Ker I_{u^\perp}$, that is, in this special
situation the kernel of $I_{u^\perp}$ is an inner ideal (namely $K_u$).
As above, $K_b$ is complemented, and hence has a natural JTS-structure.
Let us describe this structure: realize elements of $K_b$ by matrices as mentioned above;
then
$$
\begin{pmatrix} 0 & X^t \cr -X & 0 \end{pmatrix} 
\begin{pmatrix} 0 & Y^t \cr -Y & 0 \end{pmatrix} 
\begin{pmatrix} 0 & X^t \cr -X & 0 \end{pmatrix} =
\begin{pmatrix} 0 & \langle X,Y \rangle X^t \cr - \langle X,Y \rangle X & 0 \end{pmatrix} 
$$
for row vectors $X,Y \in \K^{n-1}$.
This is exactly the JTS-structure of the ``projective'' JTS $M(1,n-1;\K)$.
It follows  
also that $fgf$ is proportional to $f$ for all $f \in K_u$ and $g \in \Asym(n,\K)$
(proof: $fgf=\langle f,g \rangle f$ if $g \in K_u$, and $fgf=0$ if $g \in \Ker(K_u)$).
Thus $\K f$ is an inner ideal for all $f \in K_u$ (this property is used to define
{\em point spaces} in \cite{McC}). 
Finally, the intersection of a point space $K_u$ with a standard inner ideal $I_b$
is again a point space inner ideal; thus the inner ideals
$K_u$ can be characterized as maximal point spaces.

\msk \nin {\bf 4. Type four.}  $V^+ = V^- = \K^{p,q}$.
The  inner ideals are precisely the isotropic subspaces $I$ for the form $\beta$.
We have $\Ker I = I^\perp$, and an isotropic subspace $J$ is complementary to $I$ as
inner ideal iff  $\K^n = I \oplus J^\perp = J \oplus I^\perp$.
As JTS, the {\em proper} complementary inner ideal pairs $(I,J)$ are isomorphic to the ``projective'' JTS
$M(1,r;\K)$ with $r = \dim I \leq n/2$ (\cite{LoN94}, Example 1.14). 

\msk
Note that in all cases, with the exception of proper inner ideals for type 4 and point spaces
for type 3, inner ideals are of the same type as the ambient Jordan pair; in the two
exceptional situations, inner ideals are always isomorphic to $M(1,r;\K)$.
Finally, in low dimensions there are certain isomorphisms between various types,
see Subsection \ref{smalldim}
(in particular, the isomorphism
$\Asym(4,\R) \cong \R^{3,3}$ identifies the standard ideals in $\Asym(4,\R)$ with
the isotropic subspaces of dimension 1 and 3 in $\R^{3,3}$ and the maximal point spaces with
the isotropic subspaces of dimension 2).

\subsection{Classification of isotopes: invertible elements in structure varieties}\label{sec:Isotope}

The classification of the  {\em invertible} elements in the structure varieties of 
simple finite dimensional real Jordan pairs is known:
it is based on the complete descriptions of their structure groups given by
E.\ Neher, see \cite{Be00}, Chapter XII for tables concerning
non-exceptional spaces. 
The lists are fairly long and we will not reproduce them here.
Let us just give some comments.
The tables given in loc.\ cit.\ (especially Table XII.2.6) have the same form as the ones
to be given in the next subsection; however, they appear to be longer due to the fact
that the description of {\em invertible} homotopies is less uniform:
as a rule, instead of parametrizing by general rectangular matrices $A,B$ as
below (e.g., case 1.c), one has to take invertible matrices, which forces to
distinguish the cases $p=q$ and $p \not= q$, and instead of parametrizing by
general symmetric, Hermitian or skew-symmetric matrices (e.g., case 1.c),
one uses invertible such matrices in normal form, involving signatures in the
Hermitian case and forcing  to distinguish the cases $n$ odd and $n$ even
in the skew case.
In this sense the tables from loc.\ cit.\ form a proper subset of the tables given
below.

\subsection{Structure varieties of simple classical Jordan pairs}

\begin{thm}
For each of  the classical Jordan pairs $(V^+,V^-)$ of matrix type, the following lists
of endomorphisms $\alpha$  give complete parametrizations of
the structure variety $\Svar(V^+)$. Moreover, these parametrizations are equivariant with
respect to the natural structure group actions; therefore
a classification of homotopes $T_\alpha$ up to isomorphy can
be deduced by considering structure group orbits in the respective parameter spaces. 
\end{thm}

 \nin
{\bf
1. Jordan pairs of rectangular matrices}

\msk
\nin
{\bf 1.1} $V^+ = M(p,q;\K)$, $\K=\R,\C$:

\msk
\noindent
\begin{tabular}{lll}
label & $\alpha$ &  parameter set 
 \cr
\hline
1.a & $\alpha(X)=AXA$ & $A \in M(q,p;\K)$
\cr
1.a' & $\alpha(X)=-AXA$ & $A \in M(q,p;\K)$
\cr
1.b & $\alpha(X)=AX^tB$ & $A \in \Sym(q,\K), B \in \Sym(p,\K)$
\cr
1.c & $\alpha(X)=AX^tB$ & $A \in \Asym(q,\K), B \in \Asym(p,\K)$
\end{tabular}

\msk \nin
{\bf 1.2, case of antilinear maps:} $V^+ = M(p,q;\C)$ 

\msk
\noindent
\begin{tabular}{lll}
label & $\alpha$ &  parameter set 
 \cr
\hline
1.A & $\alpha(X)=A\overline X\overline A$ & $A \in M(q,p;\C)$
\cr
1.A' & $\alpha(X)=-A\overline X\overline A$ & $A \in M(q,p;\C)$
\cr
1.B & $\alpha(X)=A\overline X^tB$ & $A \in \Herm(q,\C), B \in \Herm(p,\C)$
\end{tabular}

\msk
\noindent
{\bf 1.3} $V^+ = M(p,q;\HHH)$:

\msk
\noindent
\begin{tabular}{lll}
label & $\alpha$ &  parameter set 
 \cr
\hline
1.3.a & $\alpha(X)=A XA$ & $A \in M(q,p;\HHH)$
\cr
1.3.a' & $\alpha(X)=-AXA$ & $A \in M(q,p;\HHH)$
\cr
1.3.b & $\alpha(X)=A\overline X^tB$ & $A \in \Herm(q,\HHH), 
B \in \Herm(p,\HHH)$
\cr
1.3.c & $\alpha(X)=A\widetilde X^tB$ & $A \in \Herm(q,\widetilde \HHH), 
B \in \Herm(p,\widetilde \HHH)$
\end{tabular}

\msk \nin
{\bf 2. Jordan pairs of symmetric matrices}

\msk
\noindent
{\bf 2.1, case of linear maps:} $V^+ = V^- = \Sym(n,\K)$, $\K = \R,\C$:

\msk
\noindent
\begin{tabular}{lll}
label & $\alpha$ &  parameter set 
 \cr
\hline
2.a & $\alpha(X)=AXA$ & $A \in \Sym(n,\K)$ 
\cr
2.a' & $\alpha(X)=-AXA$ & $A \in \Sym(n,\K)$ 
\cr
2.b & $\alpha(X)=AXA$ & $A \in \Asym(n,\K)$ 
\cr
2.b' & $\alpha(X)=-AXA$ & $A \in \Asym(n,\K)$
\end{tabular}

\msk
\noindent
{\bf 2.2, case of antilinear maps:} $V^+ = V^- = \Sym(n,\C)$

\msk
\noindent
\begin{tabular}{lll}
label & $\alpha$ &  parameter set 
 \cr
\hline
2.A & $\alpha(X)=A\overline XA$ & $A \in \Herm(n,\C))$ 
\cr
2.A' & $\alpha(X)=-A\overline XA$ & $A \in \Herm(n,\C))$ 
\end{tabular}

\msk \nin
{\bf 3. Jordan pairs of Skewsymmetric matrices}

\msk
\noindent
{\bf 3.1,  case of linear maps:} $V^+ = V^- = \Asym(n,\K)$, $\K = \R,\C$:

\msk
\noindent
\begin{tabular}{lll}
label & $\alpha$ &  parameter set 
 \cr
\hline
3.a & $\alpha(X)=AXA$ & $A \in \Asym(n,\K)$ 
\cr
3.a' & $\alpha(X)=-AXA$ & $A \in \Asym(n,\K)$ 
\cr
3.b & $\alpha(X)=AXA$ & $A \in \Sym(n,\K)$ 
\cr
3.b' & $\alpha(X)=-AXA$ & $A \in \Sym(n,\K)$
\cr
3.c & $\alpha(X)=u \otimes u^* \circ  X \circ A + A \circ X \circ u \otimes u^*$ 
& $A \in \Sym(n,\K), u \in \K^n$ 
\cr
3.c' & $\alpha(X)=-u \otimes u^* \circ  X \circ A - A \circ X \circ u \otimes u^*$ 
& $A \in \Sym(n,\K), u \in \K^n$ 
\end{tabular}

\msk
\noindent
{\bf 3.2, case of antilinear maps:} $V^+ = V^- = \Asym(n,\C)$

\msk
\noindent
\begin{tabular}{lll}
label & $\alpha$ &  parameter set 
 \cr
\hline
3.A & $\alpha(X)=A\overline XA$ & $A \in i\Herm(n,\C))$ 
\cr
3.A' & $\alpha(X)=-A\overline XA$ & $A \in i\Herm(n,\C))$ 
\cr
3.B & $\alpha(X)=u \otimes u^* \circ  \overline X \circ A + A \circ \overline X \circ u \otimes u^*$ 
& $A \in \Herm(n,\C), u \in \C^n$ 
\cr
3.B' & $\alpha(X)=-u \otimes u^* \circ  \overline X \circ A - A \circ \overline X \circ u \otimes u^*$ 
& $A \in \Herm(n,\C), u \in \C^n$ 
\end{tabular}

\msk \nin
{\bf 4. Jordan pairs of Hermitian matrices}

\msk
\noindent
{\bf 4.1}  $V^+ = V^- = \Herm(n,\C)$:

\msk
\noindent
\begin{tabular}{lll}
label & $\alpha$ &  parameter set 
 \cr
\hline
4.1.a & $\alpha(X)=AXA$ & $A \in \Herm(n,\C)$ 
\cr
4.1.a' & $\alpha(X)=-AXA$ & $A \in \Herm(n,\C)$ 
\cr
4.1.b & $\alpha(X)=A\overline X \overline A^t$ & $A \in \Sym(n,\C)$ 
\cr
4.1.b' & $\alpha(X)=-A\overline X \overline A^t$ & $A \in \Sym(n,\C)$ 
\cr
4.1.c & $\alpha(X)=A\overline X \overline A^t$ & $A \in \Asym(n,\C)$ 
\cr
4.1.c' & $\alpha(X)=-A\overline X \overline A^t$ & $A \in \Asym(n,\C)$ 
\end{tabular}

\msk
\noindent
{\bf 4.2}  $V^+ = V^- = \Herm(n,\HHH)$:

\msk
\noindent
\begin{tabular}{lll}
label & $\alpha$ &  parameter set 
 \cr
\hline
4.2.a & $\alpha(X)=AXA$ & $A \in \Herm(n,\HHH)$ 
\cr
4.2.a' & $\alpha(X)=-AXA$ & $A \in \Herm(n,\HHH)$ 
\cr
4.2.b & $\alpha(X)=A\overline X \overline A$ & $A \in \Herm(n,\HHH)$ 
\cr
4.2.b' & $\alpha(X)=-A\overline X \overline A$ & $A \in \Herm(n,\HHH)$ 
\end{tabular}

\msk
\noindent
{\bf 4.3}  $V^+ = V^- = \Herm(n,\widetilde \HHH)$:

\msk
\noindent
\begin{tabular}{lll}
label & $\alpha$ &  parameter set 
 \cr
\hline
4.3.a & $\alpha(X)=AXA$ & $A \in \Herm(n,\widetilde \HHH)$ 
\cr
4.3.a' & $\alpha(X)=-AXA$ & $A \in \Herm(n,\widetilde \HHH)$ 
\cr
4.3.b & $\alpha(X)=A\widetilde X \overline A$ & $A \in \Herm(n,\widetilde \HHH)$ 
\cr
4.3.b' & $\alpha(X)=-A\widetilde X \overline A$ & $A \in \Herm(n,\widetilde \HHH)$ 
\end{tabular}

\begin{rmk}\label{notinjrk}
The theorem says that the given parametrizations are surjective; 
 we  do not claim that they are injective:
for instance, for $\alpha$ of type 1.b the pairs of matrices
$(A,\lambda B)$ and $(\lambda A,B)$ for $\lambda \in \K$ yield the
same $\alpha$. Moreover, if $p=1$ (case of ``projective''  JTS $M(1,q;\K)$), 
Type 1.c becomes trivial (since
$\Asym(1,\K)=0$), and
Type 1.a (the map $\alpha(X)=\langle X,A \rangle A$, i.e., $\alpha = A^* \otimes A$)
becomes
a special case of Type 1.b: all structural maps are given by
$\alpha(X)=BX^t$ with symmetric $B$, that is, in the ``projective case'', structural maps are
precisely the self-adjoint operators on $\K^q$.
\end{rmk}

\nin{\it Proof of the theorem.} 
We show first that the parametrizations are well-defined, i.e., the endomorphisms $\alpha$
belong indeed to $\Svar(V^+)$.
In all cases, this is checked by elementary calculations. For convenience, let us spell out
this, e.g., for the case 1.b: let $X,Y,Z \in M(p,q;\K)$, 
and  $A \in \Sym(q,\K), B \in \Sym(p,\K)$, then
\begin{eqnarray*}
\alpha T^+(X,\alpha Y,Z)&=&
A (X AY^t B Z + Z AY^t B X)^t B
\cr
&=& (AX^tB) Y (AZ^tB) + (AZ^tB) Y (AX^tB) 
\cr
&=& T^-(\alpha X,Y,\alpha Z).
\end{eqnarray*}
For $A \in \Asym(q,\K), B \in \Asym(p,\K)$ the same
calculation applies. Equivalently, one could check that 
for any $A \in M(q,q;\K)$ and $B \in M(p,p;\K)$ the
pair of linear maps
$$
g:V^+ \to V^-, \, X \mapsto AX^t B, \, \, \quad
h:V^+ \to V^-, \, X \mapsto A^t X^t B^t 
$$
is a structural transformation from $(V^+,V^-)$ to its opposite
pair $(V^-,V^+)$, and the self-adjointness condition $g=h$ holds
if $A$ and $B$ are both symmetric or both skew-symmetric.
Similarly, for Jordan pairs of symmetric matrices, 
it is useful to observe that, for any square matrix $A$,
the map $X \mapsto AXA^t$ is a well-defined endomorphism
of $\Sym(n;\K)$ which together with
$X \mapsto A^tXA$ forms a structural pair. 
The self-adjointness condition is hence satisfied
if $A$ is symmetric or skew-symmetric.
Moreover, if $A$ is any invertible matrix, the same formula
defines an element of the structure group.
 For Hermitian matrices,
similar remarks hold with respect to maps of the form
$X \mapsto AX\overline A^t$. 
Finally, cases 3.c, 3.c', 3.B, 3.B' have a different behaviour and will be discussed 
seperately below. 

\ssk
Next let us show that the parametrizations are surjective. First of all,
we have to show that all inner ideals $I$ show
indeed up in the form $I=\im \alpha$ for suitable $\alpha$. 
For the {\em principal ideals} this is immediate: by definition, they are images of
quadratic operators $Q(A):X \mapsto AXA$, which appear in the list.
Since for Jordan algebras of symmetric or Hermitian  matrices all
inner ideals are principal, this proves our assertion in these cases.
For rectangular matrices, the inner ideal
 $I_{ab} =\{ f: F \to E \mid a \subset \ker f, \im f \subset b \}$ defined above is the image of the map
$\alpha:X \mapsto A X^t B$ where $A$ and $B$ are the symmetric matrices
describing orthogonal projection (w.r.t.\ the standard scalar product on $\K^p$ and
$\K^q$) onto $b$, resp.\ onto $a^\perp$, hence all inner ideals $I$ in
$M(p,q;\K)$ are of the form $\im \alpha$ for some $\alpha$ given in the list.
Finally, the principal inner ideals in $\Asym(n,\K)$ show up in the cases
3.a in the same way as for $\Sym(n,\K)$, whereas the point spaces show up
in case 3.c which we shall discuss below separately.

\ssk
Next, for a given inner ideal $I$ we have to show that $\Svar(I)^\times$ 
parametrizes in the way described in Lemma \ref{SvarILemma} the set of all
$\alpha$ with $\im \alpha = I$ and $\ker \alpha = \ker(I)$.
This amounts to compare, for fixed $I$, the present list with Table XII.2.6 from
\cite{Be00}.
Assume first $I = I_{ab}$ is a ``rectangular inner ideal'' in $M(p,q;\R)$.
Then we may modify the matrices $A$ and $B$ defined above by a signature (case 1.b),
or (if their size is even) by standard skew matrices (case 1.c);
this corresponds to the invertible homotopies of $I$ from Table XII.2.6, loc.\ cit.,
lines 6.1.a, 6.1.b. 
If, moreover, $I$ is ``square'' (i.e., principal), then it is also an image of homotopies from
case 1.a and 1.a', which correspond to line 1.1.b in loc.\ cit.
For $\K=\C$, complex conjugation comes in as additional invertible homotopy
(case 1.A above; Table XII.2.4 line 1.A in loc.\ cit.).
Summing up, all invertible homotopies of $I$ correspond to some $\alpha$ from
the above list.
Exactly the same pattern applies to all other cases:
the invertible homotopies of $I$ correspond in general to modify by signatures,
or by standard skew matrices in even dimension, plus adding possibly a 
complex conjugation. 
As explained in section \ref{sec:Isotope}, this corresponds exactly to the classification pattern
from \cite{Be00}, Chapter XII.  

\ssk
Let us explain now that the remaining cases 3.c, 3.c', 3.B, 3.B' describe homotopies
corresponding to point-space inner ideals in $\Asym(n,\K)$.
Let $\alpha(X)=u \otimes u^* \circ  X \circ A + A \circ X \circ u \otimes u^*$.

{\bf Claim 1:} {\em $\alpha$ is a homotopy.} 
More generally, let $\tilde \alpha(X):=AXB + BXA$ for $A,B \in \Sym(n,\K)$.
In general, $\tilde \alpha$ will not be a homotopy, but we may ask under which conditions on
$A$ and $B$ this is the case. Computing $\tilde \alpha(X)Y\tilde \alpha(X)$ and
$\tilde \alpha(X\tilde \alpha(Y)X)$, we find that equality between these two terms holds if and only if
$$
AXAYBXB + BXBYAXA = AXBYBXA  + BXAYAXB .
$$
If $B$ is a symmetric {\em rank-one} operator, i.e., an operator
$u \otimes u^*: v \mapsto \langle u,v \rangle u$, then this equality holds 
for all $X,Y \in \Asym(n,\K)$:  indeed, for all $Z \in \Asym(n,\K)$, the matrix
 $BZB$ is then skew-symmetric of rank one, 
hence must be zero, so $BXB=0$, $BYB=0$, $B(XAYAX)B=0$.  This proves Claim 1.

\ssk
{\bf Claim 2.}  {\em The image of $\alpha$ is a point space.}
Consider first the case $u=e_n$ and $A$ the diagonal matrix having last diagonal coefficient
$0$ and the others $1$. Then $u \otimes u^* = E_{nn}$, so that
$A + u \otimes u^* = 1$ is the identity matrix.
Using that $E_{nn}XE_{nn}=0$ for skew-symmetric $X$, it follows that
$$
\alpha(X)=E_{nn} X A + AX E_{nn} = X - AXA = (\id - \gamma)(X)
$$
where $\gamma(X)=AXA$ is the projection onto the principal inner ideal $I_{e_n^\perp}$
(case 3.b) with kernel $K_{e_n}$, hence $\alpha = \id - \gamma$ is the projection onto 
the maximal point space $K_{e_n}$.  
The general proof of claim 2 (not necessarily maximal point spaces) is similar.

Finally, knowing that the structure variety of $K_{e_n}$ is the one of the projective
Jordan pair, that is $\Sym(n-1,\K)$ (Remark \ref{notinjrk}), augmented by complex
conjugation in the complex case, we see that all $\alpha$ having point spaces as images
are covered by cases 3.c, 3.c', 3.B, 3.B'.
\EPf

\begin{cor}
Theorem 4.2 of Part I (\cite{BeBi}) contains a {\em complete} classification of
homotopes of classical symmetric spaces, 
in the sense that the symmetric spaces described there exhaust the list of symmetric
spaces having curvature $R_\alpha$, with $\alpha$ a homotopy of one of the classical
Jordan pairs of type 1,2 and 3; the only  exception is that in loc.\ cit.\ the symmetric spaces
corresponding to the point space inner ideals for type 3 do not appear.  
\end{cor}

\Pf
Labels in Theorem 4.2 loc.\ cit.\ are chosen such that curvature tensors $[X,Y,Z]_\alpha$ given
there correspond exactly to the LTS $R_\alpha$ corresponding to the $\alpha$-homotope
of the Jordan pair $(V^+,V^-)$. 
Thus Part I furnishes a method to calculate the standard imbeddings of these LTS,
and hence the symmetric spaces in their form $G_\alpha/H_\alpha$.
\EPf

\nin
It seems likely that the symmetric spaces belonging to the point space inner ideals 
cannot be constructed via the two-involution construction from Part I, and moreover, that they
are related to the type four series in a profound way. A better understanding of this situation 
is certainly a topic worth further work.

\subsection{Homotopes of polarized Jordan triple systems}

\begin{ex}\label{ex:Ktwo}
Let us start with a naive (and important) example:
consider the direct product of the JTS $(\K,T)$ with itself, where
$T(x,y,z)=2xyz$, whence $Q(x)y=x^2 y$. We wish to compute its structure variety $\Svar(\K \times \K)$:
write $\alpha:\K^2 \to \K^2$ as a $2 \times 2$-matrix; then 
$Q(x,y)$ is the diagonal matrix having coefficients $x^2$, $y^2$, so that the condition
$\alpha Q(X) \alpha = Q(\alpha X)$ is equivalent to
$$
\begin{pmatrix} a & b \cr c & d \end{pmatrix}
\begin{pmatrix} x^2 & 0 \cr 0 & y^2 \end{pmatrix}
\begin{pmatrix} a & b \cr c & d \end{pmatrix} =
\begin{pmatrix} (ax+by)^2 & 0 \cr 0 & (cx+dy)^2 \end{pmatrix}
$$
Taking for $X$ the two base vectors $e_1,e_2$, and assuming that $\K$ is a field, a short
computation  leads to distinguish  two cases, apart from the trivial case $\alpha = 0$:
\begin{enumerate}
\item 
$a \not= 0$; then $b=c=0$, and $\alpha$ is a diagonal matrix with $d \in \K$ arbitrary,
\item
$a = 0$ and $b\not= 0$; then $d=0$ and $c=b$, hence $\alpha$ is a multiple of the
{\em exchange map} $(x,y) \mapsto (y,x)$. 
\end{enumerate}
\nin In the first case, homotopes are simply direct products of those in $\K$; in the second
case, homotopes are {\em polarized Jordan triple systems}, corresponding 
to symmetric spaces of the kind $\Gl(2,\K)/\Gl(1,\K) \times \Gl(1,\K)$ (for $\K=\R$, this
is the one-sheeted hyperboloid, see section \ref{smalldim}).
\end{ex}

The general situation is similar as in the example:
for any Jordan pair $(V^+,V^-)$, the direct sum with its opposite Jordan
pair $(V^+ \oplus V^-,V^- \oplus V^+)$ is again a Jordan pair, and this pair 
carries a canonical involution (the exchange map), so that it becomes
a Jordan triple system $V = V^+ \oplus V^-$ with
$$
T((x,x'),(y,y'),(z,z')) = (T^+(x,y',z),T^-(x',y,z')) .
$$
Such JTS are called {\em polarized}; in the framework of $3$-graded Lie-algebras,
they correspond to the direct sum $\g\oplus \g$ with the exchange involution (where the
second copy of $\g$ carries the opposite grading of the first one).
Every structural transformation $(f:V^+ \to V^-,g:V^- \to V^+)$
from $(V^+,V^-)$ to itself  gives rise to a homotopy of the polarized JTS
$V$ via
$$
\alpha: V^+ \oplus V^- \to V^+ \oplus V^-, \quad
(u,v) \mapsto (g(v),f(u)).
$$
%More generally: the last case is the scalar extension by
%$\K[j]:=\K[X]/(X^2 - 1)$. One may look at other scalar extensions, in particular
%by dual numbers $\K[\epsilon]=\K[X]/(X^2)$.
Seen this way, the classification of all structural endomorphisms of a Jordan pair is
a subproblem of the one of classifying homotopies. 
Put differently, classifying structural transformations from Jordan pairs to themselves is
a subproblem of classifying the self-adjoint ones from Jordan pairs to their opposite pairs.
If the Jordan pair in question admits an involution, then both problems are equivalent
(since the involution can be used to reduce one to the other). 
%
%[impression: the core of Jordan theory is in ANTI-morphisms; everything that is only
%abut morphisms is more or less directly amenable to Lie theory.]
%
Note that the classification of invertible such transformations  (the structure group) is known, again by
work of E.\ Neher.
We conjecture that the following table gives indeed a complete description, 
but we will not go here into details of the proof.

\msk
\nin {\sl Table of  structural endomorphisms $(g,h)$ of a Jordan pair $(V^+,V^-)$:
 $(g:V^+ \to V^+, h:V^- \to V^-)$.}

\msk
\noindent
\begin{tabular}{lll}
label & pair of linear maps $(g,h)$ &  parameter set 
 \cr
\hline
1.a & $g(X)=AX^tB, h(X')=A^t (X')^t B^t$ & $A,B \in M(p,q;\K)$
% j'ai vŽrifiŽ. NB: si $p\not= q$, ce n'est jamais bijectif; si $p=q$, cela donne essentiellement
% les memes espaces symŽtriques que le cas 1.b
\cr
1.b & $g(X)=AXB, h(X')=BX'A$ & $A \in M(p,p;\K),B \in M(q,q;\K)$
\cr
2 & $g(X)=AXA^t,h(X')=A^t X' A$ & $A \in M(n,n;\K)$
\cr
3.a & $g(X)=AXA^t,h(X')=A^t X' A$ & $A \in M(n,n;\K)$
\cr
3.b & $\begin{matrix}
g(X)=AX \, u \otimes u^* + u \otimes u^* \,  XA^t, \cr
\quad h(X')=A^t X' \, u \otimes u^* + u \otimes u^* \, X'A
\end{matrix} $ & $A \in M(n,n;\K)$, $u \in \K^n$
\end{tabular}

\msk
\nin
For $\K=\C$ we have, as usual, $\C$-antilinear families given by the same formulas,
with a preceding complex conjugation of $X$ and $X'$.
Tables of structural endomorphisms of spaces of Hermitian matrices are similar
as for symmetric ones (type 2), replacing the transposed matrix by the conjugate transposed matrix.
The symmetric spaces corresponding to these homotopies are given in Part I,
Theorem 4.3. 

%Then: suffices to check which ones are symmetric: $g=h$

\subsection{Structure variety of an associative pair}

We can define the {\em structure variety of an associative pair},
exactly as in the Jordan case: recall from \cite{Lo75} that
an {\em associative pair} is 
given by two vector spaces $\A^\pm$ and trilinear maps
$ \A^\pm \times \A^\mp \times \A^\pm \to \A^\pm$,
$(x,y,z) \mapsto \langle xyz\rangle^\pm$ 
such that
\[
\langle xy \langle zuv\rangle^\pm \rangle^\pm =
\langle\langle xyz\rangle^\pm uv\rangle^\pm=\langle x\langle uzy\rangle^\mp v\rangle^\pm.
\]
The symmetrized product $T^\pm(x,y,z):=\langle x,y,z \rangle + \langle z,y,z \rangle$ is then a 
Jordan pair. 
If we define the {\em structure variety of the associative pair} as before by  (\ref{SvarDefEqn}),
then it will be a subset of the structure variety of the Jordan pair $T^\pm$.
However, there is a second definition of structure variety, given by the ``opposite'' condition 
\begin{equation}
\alpha \langle x,\alpha y, z\rangle = \langle \alpha z,y,\alpha x \rangle ,
\end{equation}
which gives another subset of the structure variety of the Jordan pair $T^\pm$.
In our classification, case 1 (rectangular matrices) corresponds to the associative pair
$(M(p,q;\K),M(q,p;\K))$
with product $\langle x,y,z \rangle = xyz$. Here,  cases 1.a, 1.a', 1.A, 1.A', 1.3a, 1.3.a' correspond
to the first choice of definition of the structure variety, 
and cases 1.b, 1.c, 1.B, 1.3.b, 1.3.c to the second (opposite) choice. It follows that
Theorem 3.1 classifies elements of these kinds of structure varieties.

\subsection{Jordan pairs of type four (spin factors)}

We recall the following result of Rivillis (\cite{Ri67}):

\begin{prop}\label{Type4Class}
Let $(V,T)$ be a real simple Jordan pair of Type 4 (and of dimension at least $3$)
corresponding to the bilinear form
$\beta$. Then
$\alpha$ belongs to $\Svar(V)$ if and only if
$\alpha$ is symmetric with respect to $\beta$
and $\alpha^2=\lambda \id$ with $\lambda \in \K$.
\end{prop}

\Pf (Cf.\ \cite{Ri67}.)
Note that $\alpha$ belongs to $\Svar(V)$ if and only if
$$
\beta(\alpha x,y) \alpha z +  \beta(\alpha z,y) \alpha x -
\beta(\alpha x,\alpha z)y =
\beta(x,\alpha y) \alpha z +  \beta(z,\alpha y) \alpha x -
\beta(x,z)\alpha^2 y
$$
Choosing $y$ linearly independent from $\alpha z$ and $\alpha x$, 
we see first that $\alpha^2 y$ must be proportional to $y$
and then
$\beta(\alpha x,y)=\beta(x,\alpha y)$.
This must hold for all $x,y$,  hence $\alpha$ is
symmetric. But then we get the condition
$\beta( x,\alpha^2 z)y = \beta(x,z)\alpha^2 y$ which implies that
$\alpha^2$ must be a multiple of the identity. Conversely, it
is clear that these conditions imply that $\alpha \in \Svar(V)$.
\EPf

\nin For  $\alpha^2 = 1$, we get endomorphisms that are simultaneously
symmetric
and orthogonal with respect to  $\beta$, and
for $\alpha^2=-1$, we get endomorphisms that are simultaneously symmetric
and ``anti-orthogonal'' w.r.t.\  $\beta$.
The last case can only occur if $\beta$ has signature $(n,n)$: then
it is the real (or imaginary) part of a complex quadratic form, and
 $\alpha$ is the  complex structure:
$\langle ix,iy \rangle = - \langle x,y \rangle$. 
%The corresponding
%component of the structure variety then is the symmetric space $\SO(n,n)/\SO(n,\C)$. 
The case $\alpha^2 = 0$, $\alpha \not= 0$ can appear whenever the form
$\beta$ is neither positive nor negative definite.
The symmetric spaces corresponding to $\alpha$ are described in
\cite{Ma79}. However, a global algebraic construction of these spaces, comparable 
to the ones in the other cases mentioned above, is missing.
In particular, we have the impression that the spaces corresponding to $\alpha^2 =0$
have a close relation to symmetric spaces related to point space ideals in
$\Asym(n,\R)$, see above.

%could this be sth of the kind of duality of lines and columns from the associative case?
%eg. point space somehow parametrizes the type 4 spaces, and vice versa ?

\subsection{Low dimensional cases}\label{smalldim}

In low dimension, certain isomorphisms between members of different 
families of Jordan pairs or triple systems occur; in particular, 
matrix algebras of low rank are often isomorphic to
type four algebras $\R^{p,q}$ (recall notation from the table in subsection \ref{JPclass}). 
In the following, we list such isomorphisms (following \cite{Lo75}, p. 196 ff and \cite{Be00}, Table XII.1.5) 
along with some comments
on their inner ideals and their geometry. Recall that proper inner ideals in $\R^{p,q}$ are exactly
the isotropic subspaces; thus their maximal dimension is $\min(p,q)$; in particular,
$\R^n = \R^{n,0}$ contains no non-trivial inner ideals. 
Thus homotopies in this case are either trivial or invertible, and hence every proper
contraction of an isotope is trivial (i.e., flat).
The other extreme is the case $M(1,n;\K)$: in this case, every vector subspace is
an inner ideal, and hence there are many non-trivial contractions.

Before describing the symmetric spaces, recall from Part I (\cite{BeBi}, Sections 2 and 4) 
notation and definition of homotopes of classical Lie algebras and Lie groups, as well
as the description of symmetric spaces $G_\alpha/H_\alpha$ belonging to a homotopy
$\alpha$ (\cite{BeBi}, Theorem 4.1). Recall in particular the
{\em $A$-orthogonal Lie algebra} (for an arbitrary matrix $A \in \Sym(n,\K)$)
$$
\oo_n(A; \K) := \Asym(n;\K) \, \mbox{with Lie bracket }  \,
[X,Y]_A := XAY - YAX .
$$ 
The case $n=3$ will be particularly important below: we calculate explicitly the Lie bracket
with respect to a {\em diagonal} matrix $A = dia(a,b,c)$;
let us denote this algebra by $\oo_3((a,b,c);\K)$:
$$
\Big[
\begin{pmatrix} 0 & x & y \cr -x & 0 & z \cr -y & -z & 0 \end{pmatrix},
\begin{pmatrix} 0 & x' & y' \cr -x' & 0 & z' \cr -y' & -z' & 0 \end{pmatrix} \Big]_A =
\begin{pmatrix} 0 & c(zy'-yz') & b(zx' - xz') \cr -c(zy'-yz')   & 0 & -a(xy'-x'y)\cr - b(zx' - xz') & 
a(xy'-x'y) & 0 \end{pmatrix}
$$
In other words, denoting by $e,f,k$ the natural basis in $\Asym(3;\K)$, this Lie algebra
is $\K^3$ with commutator relations
\begin{equation} \label{3Brack}
[e,f]= ck, \quad [e,k]= af, \quad [f,k]= - be
\end{equation}
For $A=dia(1,1,1)$, we get the usual orthogonal algebra
$\oo_3((1,1,1);\K)=\oo(3;\K)$; for $A=dia(-1,-1,1)$, we get $\oo_3((-1,-1,1);\K)=\oo(2,1;\K)\cong
\sL_2(\K)$; for
$A=dia(0,1,-1)$ and $A=dia(0,1,1)$, we get two solvable $3$-dimensional algebras,
$\oo_3((0,1,-1);\K)$ and $\oo_3((0,1,1);\K)$,
isomorphic to semidirect products $\K^2 \rtimes \oo(2;\K)$, resp.\
$\K^2 \rtimes \oo(1,1;\K)$.
Finally, for $A=dia(0,0,1)$, we get the $3$-dimensional Heisenberg algebra. 
Note also that, for any of these algebras, conjugation by the diagonal matrix
$I_{2,1}$ is a Lie algebra automorphism.

\subsubsection{Dimension $2$}

There are exactly three non-isomorphic two-dimensional semisimple Jordan pairs:

\ssk

\nin
{\bf A. Two-dimensional conformal geometry:} 
$V=\C=\Sym(1,\C) = M(1,1;\C) \cong  \R^2$ (isomorphic also as Jordan algebras).
The only inner ideals are $0$ and $\C$, thus homotopies are either trivial or invertible.
We have the following list of isotopes: 

$\alpha(X) = \overline X$  gives  $M_\alpha = D$ (the unit disc),

$\alpha(X)=-\overline X$ gives $M_\alpha = S^2$ (sphere, $c$-dual of $D$); 
both contract directly to flat spaces;

$\alpha(X)=X$ gives $\C^\times$ which is self $c$-dual and flat.

\ssk \nin
{\bf B. Real projective plane:} $(V^+,V^-)=(M(1,2;\R),M(2,1;\R))$.  
Every linear subspace in $V^-$ is an inner ideal, hence there are  inner ideals of dimension one, and
there exist contractions in ``two steps''.
Recall from Remark \ref{notinjrk} that in the projective case all homotopies can be put in the form 
$\alpha(X) = BX$ with a symmetric matrix $B$. Thus classification of homotopes in the
projective case amounts to classification of $\Gl(2,\R)$-orbits in $\Sym(2,\R)$:

\begin{description}
\item{B.1.} Isotopes correspond to  
the three open orbits which we represent by the three matrices
$B = 1$, $B=I_{1,1}$, $B=-1$. The corresponding isotopes are
$$
M= D, \quad M=H=\OO(2,1)/\OO(1)\times \OO(1,1) \mbox{(one-sheeted hyperboloid)},
\quad M=\R \PPP^2 =
\OO(3)/\OO(2) \times \OO(1).
$$
\item{B.2.} There are two rank-one orbits represented by diagonal matrices
$B=E_{11}=dia(1,0)$ and $B=-E_{11}= dia(-1,0)$. The corresponding homotopes are
$$
C_H = \OO_3((1,0,1))/\OO_2((1,0)) \times \OO_1(1), \quad
C_E =  \OO_3((-1,0,1))/\OO_2((-1,0)) \times \OO_1(1)
$$
where $\OO_3((a,b,c))$ is the group having Lie algebra $\oo_3((a,b,c);\R)$,
see \cite{BeBi}. We call $C_E$ the  {\em elliptic cylinder} and $C_H$
 the {\em hyperbolic cylinder} (in \cite{BDS}
the terms {\em Poincar\'e coset}, resp.\   {\em Galilei coset} had been used). 
The spaces $C_E$, resp.\ $C_H$ can be characterized as the non-flat and
non-semisimple two-dimensional symmetric spaces, $c$-dual to each other, and
topologically the former is homeomorphic to a cylinder, and the latter to a vector space.
\end{description}

\nin Contraction relations are precisely the closure relations among $\Gl(2,\R)$-orbits in 
$\Sym(2,\R)$. Point reflection at the origin corresponds to $c$-duality. 
Note that all local isomorphism classes of $2$-dimensional symmetric spaces appear in
this picture. 
In particular, the $2$-dimensional $ax+b$-group (affine groupe of the real line) is among
these spaces. Indeed,
since the  homotopy $\alpha(X)=E_{11} X = (1 \, \, 0)X(1 \, \, 0)^t$ is both of type 1.a and 1.b,
the space $C_H$ has two interpretations: it can be written as a quotient of two groups (see above), but
like all spaces of type 1.a it is also a Lie group,
seen as symmetric space, namely the group $ \Gl_{1,2}( (1 \, \, 0);\R)$.
The Lie algebra of this group, denoted by $\gL_{1,2}((1 \, 0);\R)$,  is $\R^2 = M(1,2;\R)$ with
$$
[(x,y),(x',y')] = (x,y) (1 \, 0)^t (x',y') - (x',y')(1 \, 0)^t(x,y) = (0,xy' - x'y)
$$
which is the well-known ``$ax+b$-algebra'', given by the bracket
$[e_1,e_2]=e_2$ on $\R^2$. Thus $\Gl_{1,2}((10);\R)$ is nothing but the $ax+b$-group, 
which is naturally fibered over $\R^\times$ with fiber the translation
group of $\R$.
It is quite instructive to check by a direct computation that the transvection algebra of
the $ax+b$-group is a $3$-dimensional Lie algebra, given by the bracket relations
(\ref{3Brack}) with $c=1$, $a=1$, $b=0$, hence isomorphic to the Lie algebra
$\g=\oo_3((1,0,1);\R)$, with symmetric decomposition 
$\g=\h \oplus \p$ with $\h = \R k$, $\p = \R e \oplus \R f$. 
(Note that the transvection algebra of the $ax+b$-group is strictly smaller than the
algebra of left- and right translations which is of dimension four!)
In a similar way, the $c$-dual of $C_H$ also has two different interpretations: it is
 at the same time of type 1.a' and 1.b':
$$
C_E = \Gl_{1,2}( (10);\C)/\Gl_{1,2}( (10);\R) = \OO_3((-1,0,1))/\OO_2((-1,0)) \times \OO_1(1)
$$

%
%This feature might be highly significant when comparing with \cite{BDS}:
%it seems that this group structure somehow survives in the deformed space,
%i.e., on the unit disc, in form of the solvable group acting there.

\ssk \nin
{\bf C. Polarized Geometry: direct products.}
$V= \R^{1,1}$ is not a simple Jordan algebra: it is isomorphic to $\R \times \R$.
 The structure variety of $\R \times \R$ has been computed in Example \ref{ex:Ktwo}.
Isotopes are either uninteresting (direct product of flat spaces), or isomorphic to the one-sheeted
hyperboloid $H$. Although there exist proper inner ideals (the two isotropic lines in $\R^{1,1}$),
they do not contribute to non-trivial contractions of $H$: indeed, the computations from Example
\ref{ex:Ktwo} show that homotopies with image one of these two lines lead to direct product
situations, hence to flat spaces. Summing up, no interesting contractions of $H$ arise in this case.

\subsubsection{Dimension $3$} 

There are exactly three simple Jordan pairs of dimension $3$ over $\R$:

\ssk \nin {\bf A. Riemannian conformal space:} case of the Jordan algebra $\R^3$.
Like $\R^2$, the Jordan algebra $\R^3$ has no non-trivial inner ideals, hence
there are only isotopes contracting directly to flat space. 
These isotopes are given by symmetric matrices having coefficients
$\pm 1$ on the diagonal, leading to spaces
going from the compact $S^3$ to hyperbolic space via
its indefinite isotopes. 

\msk
\nin {\bf B. Projective $3$-space:}  $M(1,3;\R) \cong \Asym(3,\R)$ (see \cite{Lo75}, p.\ 197 for this isomorphism). 
Looking at this JTS as a space of row matrices, 
the discussion is similar to the one for $M(1,2;\R)$ above: there are
 inner ideals of all dimensions, and the structure variety
  is given by all $\alpha(X)=BX$ with
   $B \in \Sym(3,\R)$; contraction relations are closure relations
of $\Gl(3,\R)$-orbits there: the
4 open orbits, corresponding the invertible matrices $B=1,I_{2,1},I_{1,2},-1$, give the isotopes
$$
\R \PPP^3 = \OO(4)/\OO(3) \times \OO(1), \quad 
\OO(2,2)/\OO(2,1) \times \OO(1), \quad
\OO(3,1)/\OO(1) \times \OO(2,1), \quad \OO(3,1)/\OO(3) \times \OO(1) .% \mbox{ (hyperbolic space). }
$$
Next we have three different orbits of rank-2 matrices and two different orbits of rank-1-matrices,
corresponding to symmetric spaces resembling the spaces $C_E$ and $C_H$ from above. 
As above, spaces associated to the two rank-1 orbits have two different interpretations: 
 semi-positive rank-one matrices correspond to homotopies
$X \mapsto AXA$ of type 1.a, hence to group cases: for $A = (1 \, 0 \, 0)$, the same computation
as above shows that this group is the affine group of $\R^2$ ($ax+b$-group).
 Its Lie algebra is $3$-dimensional with commutator relations
$[k,e]=e$, $[k,f]=f$, $[e,f]=0$, and it is not isomorphic to any of the $3$-dimensional Lie algebras
$\oo_3((a,b,c);\R)$ described above. 

Looking now at this JTS as $\Asym(3,\R)$, the picture becomes considerably more complicated:
first of all, a glance at the list of homotopes of this JTS shows that the only invertible 
homotopies are of type 3.b or 3.b' with an invertible symmetric matrix $A$.
We may choose $A=dia(a,b,c)$, so for type 3.b we get precisely 
the orthogonal groups $\OO_3((a,b,c);\R)$ introduced above.
For type 3.b', we get their $c$-duals.
This implies the four isotopes mentioned above must be isomorphic to 
$\SO(3)$, $\SO(2,1)_o$, $\SO(3,\C)/\SO(2,1)$, $\SO(3,\C)/\SO(3)$, respectively
(indeed, this could be deduced from known isomorphisms of low dimensional Lie groups).
Choosing a degenerate diagonal matrix $A=dia(a,b,c)$, we get three other group cases,
but none of them is isomorphic {\em as a group} to the $ax+b$-group mentioned above.
However, {\em as symmetric space}, the $ax+b$-group is isomorphic to $\oo_3((110);\R)$:
indeed, both the commutator relations
$[k,e]=e$, $[k,f]=f$, $[e,f]=0$ and
$[k,e]=f$, $[k,f]=e$, $[e,f]=0$ lead to the LTS relations
$[k,k,f]=f$, $[k,k,e]=e$. 
In both cases the formula for the homotopy boils down to $\alpha((x,y,z))=x$.
The case of the group $\oo_3((100);\R)$ belongs to the homotopy
$\alpha(X)=I_{11}XI_{11}$ which is zero since $X$ is skew-symmetric;
indeed, the group $\oo_3((100);\R)$ is two-step nilpotent, hence its LTS is zero. 
Note also that all inner ideals in $\Asym(3;\R)$ are  ``point spaces'' 
(since the principal inner ideal $\Asym(2;\R)$ is isomorphic to $\R$), and hence symmetric
spaces are fibered over isotopes of $\R \PPP^1$ or $\R \PPP^2$.
(For the time being, we have no good description of homotopes corresponding to point
space ideals for $\Asym(n,\R)$ in general; in particular, we do not know whether for
$n > 3$ new group cases show up.)

\msk
\nin {\bf C. Pseudo-Riemannian conformal space:} Jordan algebra  $\Sym(2,\R) \cong \R^{2,1}$.
This is the ``small Minkowski space''. Proper inner ideals are the isotropic lines.
Thus there are non-trivial contractions. This case already shows the ``generic'' feature that
the structure variety has several irreducible components.
The isotopes are:
the two conal spaces $\Gl(2,\R)/\OO(2)$, $\Gl(2,\R)/\OO(1,1)$ (type 2.a) and their $c$-duals
$\UU(2)/\OO(2)$, $\UU(1,1)/\OO(1,1)$ (type 2.a'), and the group case
$\Sp(1;\R)=\Sl(2,\R)$ (type 2.b) and its $c$-dual $\Sl(2,\C)/\Sl(2,\R)$ (type 2.b').  
Spaces of type 2.a and 2.a' have non-trivial contractions, whereas the two spaces of type 2.b, 2.b'
(belonging to another irreducible component of $\Svar(\Sym(2,\R))$) do not.

\subsubsection{Dimension $4$}

There are exactly five simple Jordan pairs of dimension four, and there are also some
interesting direct products, such as $\R^2\times \R^2 = \C \times \C = \C^2$ (which hosts the polarized space
$\Sl(2,\C)/\C^*$). For reasons of space, the following discussions of these spaces
will be sketchy and far from being exhaustive.

\ssk \nin
{\bf A.
Riemannian Conformal geometry:}
$\R^4 \cong \HHH$ (as Jordan algebras).  
Isotopes are ranking from the compact $\HHH \PPP^1 =S^4$ to hyperbolic space
and contain the group  case $\HHH^\times$; they all
contract directly to $0$. 

\ssk \nin
{\bf B. 
Conformal Artinian geometry:}
$M(2,2;\R) \cong \R^{2,2}$ (as Jordan algebras): the principal inner ideals are the
 isotropic lines (principal), the isotropic planes are inner, but not principal.
 Isotopes are:
group type space $\Gl(2,\R)$ (type 1.a), its $c$-dual $\Gl(2,\C)/\Gl(2,\R)$ (type 1.a'),
$\OO(4)/\OO(2) \times \OO(2)$,
$\OO(2,2)/\OO(2) \times \OO(2)$,
$\OO(3,1)/\OO(2) \times \OO(1,1)$,
$\OO(2,2)/\OO(1,1) \times \OO(1,1)$ (type 1.b),
$\Sp(2,\R)/\Sp(1,\R) \times \Sp(1,\R)$ (type 1.c).
The latter space has no non-trivial contractions, whereas all other have proper
contractions corresponding to $1$-dimensional inner ideals; in particular,
$\Gl(2,\R)$ contracts to groups $\Gl_2(A;\R)$.

\ssk \nin
{\bf C.
Conformal Minkowski space:}
$\Herm(2,\C)=\R^{3,1}$ (as Jordan algebras). Inner ideals are the isotropic lines; all are principal.
Isotopes are:
conal spaces $\Gl(2,\C)/\UU(2)$, $\Gl(2,\C)/\UU(1,1)$ (type 1.1.a) and their duals, the
group type spaces $\UU(2)$ and $\UU(1,1)$ (type 1.1.a'),
$\OO(2,2)/\OO(2,\C)$ (type 1.1.b; self-$c$ dual in this particular case),
$\Sp(2,\R)/\Sp(1,\C)$ (type 1.1.c), $\OO(2,2)/\Sp(1,\C)$ (type 1.1.c');
the latter two have no non-trivial contraction, whereas all other spaces have
exactly one non-trivial contraction. 

\ssk \nin {\bf D.
Projective $4$-space:}
$M(1,4;\R)$. 
The situation is similar to the case $M(1,3;\R)$ described above.
Since there is no special isomorphism in this case, this case is already ``generic'' for real
projective spaces. In particular, the only  
group cases correspond exactly to the rank-one homotopes (groups $\Gl_{1,4}(A;\R)$,
still isomorphic to an $ax+b$-group). 

\ssk \nin {\bf E.
Complex projective plane:}
$M(1,2;\C)$. 
We distinguish $\C$-linear and $\C$-antilinear homotopes:
this discussion of the $\C$-linear case follows exactly the same pattern
as the one for $M(1,2;\R)$; since $\C$ is algebraically closed, there are less
isomorphism classes of symmetric spaces ($3$ complex spaces, corresponding to
rank $B$ being $0,1$ or $2$, compared to the $6$ real spaces from above).

For $\C$-antilinear $\alpha$, let us first remark that (for $p=1$, $q=n$; here $n=2$)
types 1.A and 1.A' become special instances of type 1.B; in other words, $\C$-antilinear
homotopies are precisely the maps of the form $F(z)=f(\overline z)$ with Hermitian
$f:\C^2 \to \C^2$. 
For invertible $f$, we get according to the various possible signatures, the isotopes
$$
\UU(3)/\UU(2) \times \UU(1), \quad
\UU(2,1)/\UU(2) \times \UU(1), \quad
\UU(2,1)/\UU(1) \times \UU(1,1) .
$$
If $f$ is of rank one, then the spaces have two interpretations:
type 1.A corresponds to $f$ of rank one and non-negative, and
type 1.A' to $f$ of rank one and non-positive. For the corresponding symmetric spaces
we have for type 1.A:
$$
\UU_3((1 , -1 ,0)/\UU_1(1) \times \UU_2((1, 0)) \cong
\Gl_{2,1}( (1 ,0); M(2,2;\R)) / \Gl_{2,1}( (1, 0;\C)\;.
$$

\subsubsection{Dimension $6$}

The cases of dimension 6 and 10 are of particular interest, since
 low dimensional isomorphisms appear there.
We just give the list of simple JTS and of special isomorphisms in dimension 6, and leave a more
detailed discussion of the corresponding homotopes and their specific features for subsequent work. 

\ssk \nin
{\bf A. Riemannian conformal space:} $\R^6$

\ssk \nin
{\bf B. Lorentzian conformal space:} $\R^{5,1} \cong \Herm(2,\H)$

\ssk \nin
{\bf C. Artinian conformal space:}
$\Asym(4,\R) = \R^{3,3}$ (see \cite{Lo75}, p. 200, for this isomorphism)
% for an explicit construction of this isomorphism:
%the quadratic form comes the Pfaffian)
%
% standard inner ideals: isotropic lines and $3$-spaces;
%point spaces: isotropic planes.
%Among the isotopes are the group $\SO(4)$ (not simple as group!), $\SO(3,1)$,
%$\SO(2,2)$, they have (as groups!) a very rich contraction structure, corresponding
%to orbits of $\Gl(4,\R)$ in $\Sym(4,\R)$, up to sign. 
%Generally, orthogonal groups have a richer contraction structure than symplectic ones
%(especially in small dimension) since their parameter space is given by symmetric
%matrices instead of skew-symmetric matrices. 

\ssk \nin
{\bf D. Pseudo-Riemannian conformal space of signature $(4,2)$:} $\R^{2,4}$ % \cong ?

\ssk \nin
{\bf E. Complex projective space:}
$M(1,3;\C) \cong \Asym(3,\C)$

\ssk \nin
{\bf F. Complex Lagrangian geometry of a symplectic form on $\C^4$:}
 $\Sym(2,\C) \cong \C^3$ 

\ssk \nin
{\bf G. Real Lagrangian geometry of a symplectic form on $\R^6$:}
$\Sym(3,\R)$

\ssk \nin {\bf H. Real Grassmannian of $2$-spaces in $\R^5$:}
$M(2,3;\R)$
% inner ideals of dimension 1,2,3,4;
%groupe cases correspond to dimension $1$ and $4$ (?)

%{\bf Dimension $10$:}

%check whether $\Herm(2,\tilde \HHH)= \R^{5,5}$ !?

%recall $\Herm(2,Octon)=\R^{1,9}$ 

%

\end{document}